\documentclass[%
 aip,
 amsmath,amssymb,
preprint,%
]{revtex4-1}

\usepackage{graphicx}
\usepackage{dcolumn}
\usepackage{bm}
\usepackage{textcomp}
\usepackage{etoolbox}
\usepackage{amsmath}
\usepackage{float}

\usepackage{subfigure}
\usepackage{mathrsfs} 

\usepackage{amsthm}

\def\Dt{\De t}
\def\Dx{\De x}
\def\Dy{\De y}

\def\bn{\mathbf{n}}

\def\bu{\mathbf{u}}

\def\bx{\mathbf{x}}

\def\bA{\mathbf{A}}
\def\bB{\mathbf{B}}

\def\bF{\mathbf{F}}
\def\bG{\mathbf{G}}

\def\bI{\mathbf{I}}

\def\bK{\mathbf{K}}
\def\bL{\mathbf{L}}

\def\bR{\mathbf{R}}

\def\bU{\mathbf{U}}
\def\bV{\mathbf{V}}
\def\bW{\mathbf{W}}

\def\n{\noindent}
\def\pt{\partial}

\def\minmod{\mbox{minmod}}
\def\f#1#2{\frac {#1}{#2}}
\def\f32{\frac 32}

\def\d{\displaystyle}

\def\beq{\begin{equation}}
\def\eeq{\end{equation}}
\def\bga{\begin{array}}
\def\eda{\end{array}}

\def\De{\Delta}

\def\gm{\gamma}

\def\al{\alpha}

\def\Om{\Omega}

\def\d{\displaystyle}
\def\dfr#1#2{\displaystyle{\frac{#1}{#2}}}

\def\beql#1{\begin{equation}\label{#1}}
\def\beqs{\begin{equation*}}
\def\eeqs{\end{equation*}}

 \newtheorem{thm}{Theorem}[section]
 \theoremstyle{remark}
 \newtheorem{rem}[thm]{Remark}

\makeatletter
\def\@email#1#2{%
 \endgroup
 \patchcmd{\titleblock@produce}
  {\frontmatter@RRAPformat}
  {\frontmatter@RRAPformat{\produce@RRAP{*#1\href{mailto:#2}{#2}}}\frontmatter@RRAPformat}
  {}{}
}%

\makeatother
\begin{document}

\renewcommand{\tablename}{TABLE}
\renewcommand{\figurename}{FIG}

\preprint{\emph{Preprint submitted to Physics of Fluids}}
\title[]{Generalized Riemann Problem Method for the Kapila Model of Compressible Multiphase Flows}

\author{Tuowei Chen}
\affiliation{Department of Mathematics, Hong Kong University of Science and Technology, Clear Water Bay, Kowloon, HongKong}%
\author{Zhifang Du*}%
 \email{du\_zhifang@iapcm.ac.cn}
 \affiliation{Institute of Applied Physics and Computational Mathematics, \\ 100048 Beijing, PR China}%


\begin{abstract}
A second-order accurate and robust numerical scheme is developed for the Kapila model to simulate compressible multiphase flows. The scheme is formulated within the finite volume framework with the generalized Riemann problem (GRP) solver employed as the cornerstone. Besides Riemann solutions, the GRP solver provides time derivatives of flow variables at cell interfaces, achieving second-order accuracy in time within a single stage. The use of the GRP solver enhances the capability of the resulting scheme to handle the stiffness of the Kapila model in two ways. First, the coupled values, i.e., Riemann solutions and time derivatives, give the cell interface values of flow variables at the new time level, yielding an approximation to the cell average of the velocity divergence at the new time level in a computational step. This allows a semi-implicit time discretization to the stiff source term of the volume fraction equation. Second, the effects of source terms are directly included in the numerical flux via the computation of time derivatives. The resulting numerical flux is able to capture the physics of interactions between phases, and the robustness of the scheme is therefore further improved. Several challenging numerical experiments are conducted to demonstrate the good performance of the proposed finite volume scheme. In particular, a test case with a nonlinear smooth solution is designed to verify the numerical accuracy.
\end{abstract}

\maketitle

\section{Introduction}\label{sec:intro}
Compressible multiphase flows appear in numerous scientific and engineering disciplines, including astrophysics \cite{kamaya-1996,star-mre,combustion},
cavitation flows \cite{cavitation-saurel-jcp,saurel-jfm,saurel-ijmf},
and deflagration-to-detonation transitions \cite{detonation-saurel-jfm,detonation-saurel-jcp,detonation-saurel-sw}. 
There has been an enormous amount of successful applications of sharp interface methods in numerical simulations involving compressible multiphase flows
\cite{Tryggvason-2009,LiXL-2015,Fedkiw,Osher-2018}.
In contrast, this work concentrates on the development of a shock-capturing scheme that is capable of capturing interfaces between pure fluids as well as fluid mixtures.
To achieve this, the Kapila model is employed in the present work for its capability of effectively modeling the underlying physical process of compressible multiphase flows \cite{kapila-structure,kapila,review-saurel}.

In Ref.~\onlinecite{b-n}, Baer and Nunziato proposed a seven-equation model, known as the Baer-Nunziato (BN) model, to describe the dynamics of multiphase flows, specifically in solid granular explosives mixed with gaseous products.
The BN model includes mass, momentum, and energy balance laws for each individual phase. 
The six governing equations are augmented by a seventh one modeling the gaseous volume fraction evolution. 
After being further extended to modeling general multiphase flows \cite{bn-abgrall}, 
much attention has been devoted to developing numerical schemes for the seven-equation model \cite{bn-saurel-warnecke,bn-warnecke,bn-kapila,bn-karni,bn-lei}.
A five-equation model, known as the Kapila model, was proposed in Ref.~\onlinecite{kapila} as a reduced version of the well-known BN model. This reduction is based on the assumption that, in many practical applications, the interphase exchanges of momentum and kinetic energy occur rapidly. The Kapila model emerges as the asymptotic limit in cases of strong interphase drag and negligible compaction viscosity \cite{kapila}.
The Kapila model consists of continuity equations of each phase, conservation laws of the total momentum and the total energy, and the volume fraction equation with a stiff source term that represents the interaction between phases. The sound speed of the mixture is defined by the Wood formula \cite{wood}.
The same five-equation model was also derived in Ref.~\onlinecite{kapila-puckett}.
Although this model was derived from the seven-equation model for mixtures, it can be applied to numerical simulations of interfaces separating either pure fluids or fluid mixtures \cite{review-saurel}.

After its introduction, the Kapila model has been widely used in numerical simulations of compressible multiphase flows \cite{murrone,kapila-netherland,kapila-tian}.
Nevertheless, such efforts face two distinct difficulties.
In the first place, the stiff source term in the volume fraction equation, which arises as the result of the stiff limit of the relaxation process in the original seven-equation model, is particularly difficult to approximate and makes the bound-preserving property of the volume fraction a challenging task.
Furthermore, the non-monotonic behavior of the Wood sound speed in relation to volume fractions leads to Mach number oscillations in the numerical diffusion zone of an interface. 
To mitigate these difficulties, Saurel et al. proposed a projection-relaxation scheme \cite{saurel-projection-I, saurel-projection-II}, which relies on resolving subcell structures of the numerical solution.
Later, this method was further enhanced by relaxing the pressure equilibrium assumption \cite{saurel-simple}.
 
The present paper provides an alternative approach to address the stiffness of the Kapila model.
To accomplish this, we highlight the significance of the generalized Riemann problem (GRP) solver in developing high-resolution and robust numerical schemes for compressible fluid flows. 
The GRP solver is a typical representative of Lax-Wendroff type flux solvers, which allows the corresponding numerical scheme to achieve second-order accuracy in time within a single-stage, comparing with the Runge-Kutta type schemes.
The GRP solver was initially proposed in the context of gas dynamics \cite{benartzi-84,benartzi-01},
and a more straightforward and efficient version of the GRP solver was later developed in the Eulerian framework \cite{grp-06}.
Under the acoustic assumption, the GRP can be linearized and thus solved approximately \cite{menshov-1991,Toro-ADER,Toro-ADER-pnpm,Toro-ADER-blood}.
Developing a GRP solver for the Kapila model to address its stiffness is the focus of this paper.


In comparison with the Riemann problem, the GRP is formulated as an initial value problem (IVP) of hyperbolic partial differential equations (PDEs) with piecewise smooth initial data in place of piecewise constant ones used in the Riemann problem. 
The GRP features waves emanating from the initial discontinuity that no longer travel along straight lines, and provides instantaneous time derivatives of flow variables along with Riemann solutions.
The fundamental approach adopted to solve the GRP involves the so-called Lax-Wendroff procedure \cite{L-W}, which employs governing equations to quantify the relation between the time derivative of the solution and its initial spatial variation \cite{tsc-li}. 
The coupled values, i.e., Riemann solutions and instantaneous time derivatives,
provide evolutionary cell interface values of flow variables.
This feature has been successfully used to develop compact data reconstructions \cite{du-hweno,gks-hweno-0,gks-hweno}, which rely on approximating gradients of flow variables at the new time level in a computational step.
Additionally, it is illustrated that time derivatives incorporate the transversal effect of the flow field into numerical fluxes, which facilitates an accurate resolution to multi-dimensional effects in numerical computations  \cite{Lei-Li-2018-AMM}. 

Regarding the Kapila model, the application of the GRP solver enhances the robustness of the numerical scheme in two ways.
The first notable advantage is that the source term in the volume fraction equation can be evolved using a semi-implicit discretization.
This is possible since the GRP solver provides cell interface values of the flow variables at the new time level, which, together with the Gauss-Green formula, give the cell average of velocity divergence at the new time level.
When specifically implementing the Crank-Nicolson evolution,
a time step restriction to preserve the bound of volume fractions can be established.
Second, thanks to the Lax-Wendroff procedure, which makes thorough use of governing PDEs, the compaction or expansion effect carried by the source term in the volume fraction equation is fully resolved by utilizing time derivatives at cell interfaces.
As a result, the capability of numerical fluxes to capture interactions between phases is considerably improved.

This paper is organized as follows.
The GRP-based finite volume scheme for the Kapila model is constructed in Section \ref{sec:scheme}.
In addition, the time step restriction to preserve the bound of volume fractions is established therein. 
The corresponding GRP solver is developed in Section \ref{sec:grp}.
Finally, we present results of several numerical experiments in Section \ref{sec:numer} to demonstrate the good performance of the proposed scheme. Several concluding remarks are put in Section \ref{sec:discussions}.

\vspace{2mm}
\section{Finite Volume Scheme for the Kapila Model}\label{sec:scheme}
This section is dedicated to the development of a finite volume scheme for the Kapila model.
In two-phase cases, the governing equations of the Kapila model are
\begin{subequations}\label{eq:govern}
\begin{align}
&\dfr{\pt\zeta_1\rho}{\pt t}+\nabla\cdot(\zeta_1\rho\bu)=0,    \label{eq:zeta}\\
&\dfr{\pt\rho}{\pt t}+\nabla\cdot(\rho\bu)=0,         \label{eq:rho}\\
&\dfr{\pt\rho\bu}{\pt t}+\nabla\cdot(\rho\bu\otimes\bu)+\nabla p=0,    \label{eq:mom}\\
&\dfr{\pt\rho E}{\pt t}+\nabla\cdot[(\rho E+p)\bu]=0,        \label{eq:ene}\\
&\dfr{\pt\al_1}{\pt t}+\bu\cdot\nabla\al_1=
-\al_1\al_2
\dfr{\rho_1{c_1}^2-\rho_2{c_2}^2}{\al_1\rho_2{c_2}^2+\al_2\rho_1{c_1}^2}
\nabla\cdot\bu,    \label{eq:vof}
\end{align}
\end{subequations}
where $\rho$ is the mixture density, $E$ is the total energy per unit mass, $p$ is the pressure, and $\bu$ is the velocity.
In two space dimensions, $\bu=(u,v)^\top$.
For each phase $k$ $(k=1,2)$, $\al_k$ and $\zeta_k$ represent the volume and mass fractions, satisfying
\begin{equation}\label{eq:saturation}
\al_1+\al_2=1, \ \ \zeta_1+\zeta_2=1.
\end{equation}
Note that given the continuity equation \eqref{eq:zeta} for the phase $k=1$, the other one for the phase $k=2$ is replaced by the conservation law of the mixture mass \eqref{eq:rho}.
The phase density $\rho_k$ is defined by
\beqs
\rho_k=\dfr{\zeta_k}{\al_k}\rho,
\eeqs
which leads to $\rho=\al_1\rho_1+\al_2\rho_2$.
The total energy of the mixture is
\beql{eq:total-ene-dfr}
E=\dfr 12 \left|\bu\right|^2+e,
\eeq
where $e$ is the specific internal energy of the mixture. 
The thermodynamical variables of each phase follow its Gibbs relation
\beql{eq:gibbs}
de_k=T_kds_k-pd\tau_k,
\eeq
where $e_k$ and $T_k$ are the specific internal energy and the temperature of each phase $k$ and  $\tau_k=\frac{1}{\rho_k}$ is the specific volume.
From Gibbs relation \eqref{eq:gibbs}, the specific internal energy of each phase depends on other thermodynamical variables via the equation of state (EOS)
\beql{eq:eos-e}
e_k=e_k(p,\rho_k).
\eeq
In this paper, we consider the stiffened gas EOS:
\beql{eq:eos-sg}
p=(\gm_k-1)\rho_ke_k-\gm_k\pi_k,
\eeq
where $\gm_k>1$ is the specific heat ratio and $\pi_k\geq0$ is a material parameter. 
The sound speed of each phase $k$ is defined as
\beqs
{c_k}^2=\left(\dfr{\pt p}{\pt\rho_k}\right)_{s_k},
\eeqs
where $s_k$ is the entropy.
The mixture sound speed is defined by the Wood formula \cite{wood} as
\beql{eq:cs-wood}
\dfr{1}{\rho c^2}=\dfr{\al_1}{\rho_1 {c_1}^2}+\dfr{\al_2}{\rho_2 {c_2}^2}.
\eeq
For stiffened gases, the sound speed of each phase $k$ is explicitly defined by
\beql{eq:cs-def-stiffened}
{c_k}^2=\dfr{\gm_k(p+\pi_k)}{\rho_k}.
\eeq

In the remainder of this section, the GRP-based finite volume scheme for the Kapila model together with the corresponding semi-implicit time discretization of the volume fraction equation is developed. Afterwards, this section is closed by analyzing the bound-preserving property of volume fractions.
To make the present paper self-containing, the linear data reconstruction procedure employed in Refs.~\onlinecite{benartzi-01,grp-06} is put into Appendix \ref{sec:app}.

\vspace{2mm}
\subsection{Finite volume discretization of governing equations} \label{subsec:fv}
In order to apply the finite volume scheme for hyperbolic balance laws, add $\al_1\nabla\cdot\bu$ on both sides of \eqref{eq:vof} and rewrite it as
\beq\label{eq:vof-balance}
\bga{l}
\dfr{\pt\al_1}{\pt t}+\nabla\cdot(\al_1\bu)=
\dfr{\al_1\rho_2{c_2}^2}{\al_1\rho_2{c_2}^2+(1-\al_1)\rho_1{c_1}^2}\nabla\cdot\bu.
\eda
\eeq
With \eqref{eq:vof-balance}, rewrite the governing equations \eqref{eq:zeta}-\eqref{eq:vof}  in two space dimensions as 
\beql{eq:govern-balance}
\dfr{\pt\bU}{\pt t}+\dfr{\pt\bF(\bU)}{\pt x}+\dfr{\pt\bG(\bU)}{\pt y}=\mathbf{K}(\bU,\nabla\bU),
\eeq
where
\beql{eq:govern-def}
\bga{l}
\bU=[\zeta_1\rho,\rho,\rho u, \rho v, \rho E, \al_1]^\top,\\[2mm]
\bF=[\zeta_1\rho u,\rho u,\rho u^2+p, \rho uv, (\rho E+p)u, u\al_1]^\top,\\[2mm]
\bG=[\zeta_1\rho v,\rho v,\rho uv, \rho v^2+p, (\rho E+p)v, v\al_1]^\top,\\[2mm]
\mathbf{K}=[0,0,0,0,0,\al_1K\left(\dfr{\pt u}{\pt x}+\dfr{\pt v}{\pt y}\right)]^\top,
\eda
\eeq
and
\beql{eq:Xi-def}
K=\dfr{\rho_2{c_2}^2}{\al_1\rho_2{c_2}^2+(1-\al_1)\rho_1{c_1}^2}.
\eeq
Divide the computational domain $\Om$ into uniform rectangular cells $\{\Om_{ij}\}$ as
\beqs
\Om=\bigcup_{i,j}\Om_{ij},
\eeqs
where
\beqs
\Om_{ij}=[x_{i-\frac 12},x_{i+\frac 12}]\times[y_{j-\frac 12},y_{j+\frac 12}].
\eeqs
The cell size is denoted as $\Dx=x_{i+\frac 12}-x_{i-\frac 12}$ and $\Dy=y_{j+\frac 12}-y_{j-\frac 12}$. Denote $x_i=\frac 12(x_{i-\frac 12}+x_{i+\frac 12})$ and  $y_j=\frac 12(y_{j-\frac 12}+y_{j+\frac 12})$.
By integrating governing equations \eqref{eq:govern-balance} in the space-time control volume $\Omega_{ij}\times[t^n,t^{n+1}]$, the GRP-based finite volume scheme for \eqref{eq:govern-balance} is
\beql{eq:scheme}
\bga{l}
\overline\bU^{n+1}_{ij}=\overline\bU^{n}_{ij}
-\dfr{\Dt}{\Dx}\left(\bF^{n+\frac 12}_{i+\frac 12,j}-\bF^{n+\frac 12}_{i-\frac 12,j}\right)
-\dfr{\Dt}{\Dy}\left(\bG^{n+\frac 12}_{i,j+\frac 12}-\bG^{n+\frac 12}_{i,j-\frac 12}\right)\\[3mm]
\qquad\qquad\ \ \ \ \
+\Dt\left[
(1-C_\text{im})\mathbf{K}^n_{ij}
+C_\text{im}\mathbf{K}^{n+1}_{ij}
\right],
\eda
\eeq
where $\Delta t=t^{n+1}-t^n$ is the time step and
\beqs
\overline{\bU}^n_{ij}=\dfr{1}{|\Omega_{ij}|}\int_{\Omega_{ij}}\bU(\bx,t^n)d\bx, \ \ \
\overline{\bU}^{n+1}_{ij}=\dfr{1}{|\Omega_{ij}|}\int_{\Omega_{ij}}\bU(\bx,t^{n+1})d\bx.
\eeqs
The finite volume scheme \eqref{eq:scheme} distinguishes from those developed solely based on Riemann solutions by two means: numerical fluxes are  evaluated at $t^{n+\frac 12}=t^n+\frac 12\Delta t$ and a semi-implicit time discretization to the source term is used.

The numerical flux across the cell interface $\Gamma_{i+\frac 12,j}=\{x_{i+\frac 12}\}\times[y_{j-\frac 12},y_{j+\frac 12}]$ is approximated with second-order accuracy by
\beql{eq:flux-def}
\bF^{n+\frac 12}_{i+\frac 12,j}=\bF(\bU(\bx_{i+\frac 12,j},t^{n+\frac 12})),
\eeq
where $\bx_{i+\frac 12,j}=(x_{i+\frac 12},y_{j})$.
The cell interface value $\bU(\bx_{i+\frac 12,j},\cdot)$ at $t^{n+\frac 12}$ is obtained by the Taylor expansion
\beq\label{eq:mid-time}
\bU(\bx_{i+\frac 12,j},t^{n+\frac 12})
=
\bU^{n,*}_{i+\frac 12,j}+\dfr{\Dt}{2}\left(\dfr{\pt\bU}{\pt t}\right)^{n,*}_{i+\frac 12,j},
\eeq
which is second-order accurate in time. The Riemann solution $\bU^{n,*}_{i+\frac 12,j}$ and the instantaneous time derivative $(\pt\bU/\pt t)^{n,*}_{i+\frac 12,j}$ are obtained by solving the GRP of governing equations \eqref{eq:govern-balance} at $(\bx_{i+\frac 12,j},t^n)$, which is formulated as the IVP
\beql{eq:grp-def-local}
\bga{l}
\dfr{\pt\bU}{\pt t}+\dfr{\pt\mathbf{F}(\bU)}{\pt x}+\dfr{\pt\mathbf{G}(\bU)}{\pt y}=\bK(\bU,\nabla\bU),\\[2mm]
\bU(\bx,t^n)=\left\{
\bga{ll}
\bU^n_{i,j}(\bx), & x<x_{i+\frac 12},\\
\bU^n_{i+1,j}(\bx), & x>x_{i+\frac 12},
\eda
\right.
\eda
\eeq
where $\bx=(x,y)$.
The piecewise smooth initial data $\bU^n_{i,j}(\bx)$ and $\bU^n_{i+1,j}(\bx)$ are obtained by the space data reconstruction procedure presented in Appendix \ref{sec:app}.
The GRP solver to solve the IVP \eqref{eq:grp-def-local} will be developed in the next section.
If a higher order of the spatial accuracy is desired, the numerical fluxes can be computed by the Gauss-Legendre integration and we do not elaborate tedious but basic computations.
The numerical flux $\bG_{i,j+\frac 12}^{n+\frac 12}$ can be computed in the same manner.

The parameter $C_\text{im}\in[0,1]$ determines the time discretization to the source term in the volume fraction discretization. The source term $\bK$ is approximated at $t^n$ and $t^{n+1}$ by
\beqs
\mathbf{K}_{ij}^n=[0, 0, 0, 0, 0,\overline{(\al_1)}_{ij}^n K_{ij}^n\eta_{ij}^n]^\top, \ \ 
\mathbf{K}_{ij}^{n+1}=[0, 0, 0, 0, 0,\overline{(\al_1)}_{ij}^{n+1} K_{ij}^{n+1}\eta_{ij}^{n+1}]^\top,
\eeqs
where
\beq\label{eq:div}
\eta_{ij}^n=\dfr{1}{|\Omega_{ij}|}\int_{\Omega_{ij}}\left.\nabla\cdot\bu\right|_{t=t^n}d\bx, \ \ \
\eta_{ij}^{n+1}=\dfr{1}{|\Omega_{ij}|}\int_{\Omega_{ij}}\left.\nabla\cdot\bu\right|_{t=t^{n+1}}d\bx.
\eeq
The velocity divergence in \eqref{eq:div} is approximated by the Gauss-Green formula, and details will be specified in the next subsection.

\vspace{2mm}
\subsection{Crank-Nicolson time evolution of the volume fraction equation} \label{subsec:cn-source}
For the simplicity of the presentation, this subsection illustrates the Crank-Nicolson evolution to the volume fraction equation \eqref{eq:vof-balance} in one space dimension.
In order to further simplify notations, we drop the subscript referring to the phase index and denote the volume fraction of the phase $k=1$ by $\al$.

For arbitrary $C_\text{im}\in[0,1]$, the semi-implicit discretization to \eqref{eq:vof-balance} in the cell $I_i=[x_{i-\frac 12}, x_{i+\frac 12}]$  is
\beql{eq:vof-balance-CN}
\bga{l}
\bar{\al}_{i}^{n+1}=\bar{\al}_{i}^{n}-\dfr{\Dt}{\Dx}\left[(u\al)^{n+\frac 12}_{i+\frac 12}-(u\al)^{n+\frac 12}_{i-\frac 12}\right]
\\[3mm]
\qquad\qquad\quad
+(1-C_\text{im})\Dt
\bar{\al}_{i}^{n}K^n_{i}
\eta_{i}^n
+C_\text{im}\Dt
\bar{\al}_{i}^{n+1}K^{n+1}_{i}
\eta_{i}^{n+1}.
\eda
\eeq
By using the Remann solution of the velocity $u_{i+\frac 12}^{n,*}$ and the Newton-Leibniz formula, the approximation to the divergence of the velocity field at $t^n$ is 
\beql{eq:divergence-expl}
\eta_i^n=\dfr{u_{i+\frac 12}^{n,*}-u_{i-\frac 12}^{n,*}}{\Dx}.
\eeq
Similarly,
\beql{eq:divergence-impl}
\eta_{i}^{n+1}=\dfr{\hat{u}_{i+\frac 12}^{n+1}-\hat{u}_{i-\frac 12}^{n+1}}{\Dx}, 
\eeq
where $\hat{u}_{i+\frac 12}^{n+1}=u_{i+\frac 12}^{n,*}+\Dt(\pt u/\pt t)_{i+\frac 12}^{n,*}$ is the estimation of $u$ at $(x_{i+\frac 12},t^{n+1})$.

\vspace{2mm}
\begin{rem}
The time derivative  $(\pt u/\pt t)^{n,*}_{i+\frac 12}$ allows one to evolve the solution at cell interfaces. As a result, the velocity divergence at $t^{n+1}$ can be approximated
in a straightforward way \eqref{eq:divergence-impl}. This is not available for finite volume schemes solely based on Riemann solutions.
\end{rem}

\vspace{2mm}
\begin{rem}
In two space dimensions, the Gauss-Green formula is used. Therefore, the approximation to the divergence of the velocity field in a polygonal cell $\Omega_J$ at $t^n$ is
\beq
\eta_J^{n}=\dfr{1}{|\Omega_J|}\sum_{\{\Gamma\}}|\Gamma|
\Big(\sum_g^Gw_g\bu^{n,*}_{\Gamma,g}\cdot\bn_{\Gamma}\Big),
\eeq
where $\{\Gamma\}$ is the set of all edges of $\Omega_J$, $\bn_{\Gamma}$ is the unit outer normal of $\Gamma$, $\bu^{n,*}_{\Gamma,g}$ is the Riemann solution of $\bu$ at the $g$-th Gauss quadrature point on $\Gamma$, and $w_g$ is the corresponding quadrature weight. In particular, given a rectangular cell $\Omega_{ij}$, we obtain the second-order accurate approximation:
\beq
\eta_{ij}^{n}=
\dfr{u^{n,*}_{i+\frac 12,j}-u^{n,*}_{i-\frac 12,j}}{\Delta x}
+
\dfr{v^{n,*}_{i,j+\frac 12}-v^{n,*}_{i,j-\frac 12}}{\Delta y}
.
\eeq
The divergence approximation $\eta_{ij}^{n+1}$ can be obtained in the same way by using $\hat{u}^{n+1}_{i\pm\frac 12,j}$ and $\hat{v}^{n+1}_{i,j\pm\frac 12}$.
\end{rem}

By taking $C_\text{im}=\frac{1}{2}$, the Crank-Nicolson discretization of the volume fraction equation is obtained.
In such a case, rewrite \eqref{eq:vof-balance-CN} as
\beql{eq:vof-balance-CN-2}
\bar{\al}_{i}^{n+1}
=\tilde{\al}_{i}^{n}+\theta
\bar{\al}_{i}^{n+1}K^{n+1}_{i},
\eeq
where $\theta = \frac{\Dt}{2}\eta_{i}^{n+1}$ and
\beql{eq:def-alpha-star}
\bga{l}
\tilde{\al}_{i}^{n}\triangleq\bar{\al}_{i}^{n}-\dfr{\Dt}{\Dx}\left[(u\al)^{n+\frac 12}_{i+\frac 12}-(u\al)^{n+\frac 12}_{i-\frac 12}\right]
+\dfr{\Dt}{2}
\bar{\al}_{i}^{n}
K^n_{i}
\eta_{i}^n.
\eda
\eeq
In the case of stiffened gas EOS, \eqref{eq:Xi-def} yields
\beq\label{eq:Omega-def}
K^{n+1}_{i}=\dfr{\gm_2(p^{n+1}_{i}+\pi_2)}
{\bar{\al}^{n+1}_{i}\gm_2(p^{n+1}_{i}+\pi_2) 
+
\left[1-\bar{\al}^{n+1}_{i}\right]\gm_1(p^{n+1}_{i}+\pi_1)}.
\eeq
Here the pressure $p_i^{n+1}$ is determined by the EOS of the mixture, which depends on $\bar\al_i^{n+1}$.
Given the internal energy $\overline{(\rho e)}^{n+1}_i$, which is obtained
by the finite volume scheme \eqref{eq:scheme}, the pressure is determined by the relation
\beqs
\bga{l}
\overline{(\rho e)}^{n+1}_i=\bar{\al}^{n+1}_{i}\overline{(\rho_1 e_1)}^{n+1}_i+(1-\bar{\al}^{n+1}_{i})\overline{(\rho_2 e_2)}^{n+1}_i\\[2mm]
\quad\quad\ \ \ \ \ =\dfr{\bar{\al}^{n+1}_{i}(p_i^{n+1}+\gm_1\pi_1)}{\gm_1-1}
+\dfr{(1-\bar{\al}^{n+1}_{i})(p_i^{n+1}+\gm_2\pi_2)}{\gm_2-1}.
\eda
\eeqs
Therefore, the pressure $p_i^{n+1}$ depends on the volume fraction $\bar\al_i^{n+1}$ through
\beq\label{eq:pres-impl}
\bga{l}
p^{n+1}_i=
\dfr{(\gm_1-1)(\gm_2-1)\overline{(\rho e)}^{n+1}_i
-\gm_1\pi_1(\gm_2-1)\bar{\al}^{n+1}_{i}
-\gm_2\pi_2(\gm_1-1)\left[1-\bar{\al}^{n+1}_{i}\right]}
{(\gm_2-\gm_1)\bar{\al}^{n+1}_{i} + \gm_1-1}\\[3mm]
\qquad\  
\triangleq
\dfr{L_1\left(\bar{\al}^{n+1}_{i}\right)}{L_2\left(\bar{\al}^{n+1}_{i}\right)}.
\eda
\eeq
Substituting   \eqref{eq:def-alpha-star}, \eqref{eq:Omega-def}, and \eqref{eq:pres-impl} into \eqref{eq:vof-balance-CN-2} derives a semi-implicit evolution of the volume fraction equation \eqref{eq:vof-balance}, which is written as the nonlinear equation
\beq\label{eq:f}
f\left(\bar{\al}^{n+1}_{j};\tilde{\al}_{j}^{n}\right)=0.
\eeq
The left hand side is a rational function defined as
\beq\label{eq:def-f}
f(\xi;\tilde\al)=\xi-\tilde\al- \dfr{\theta\big[\gm_2L_1(\xi)+\gm_2\pi_2L_2(\xi)\big]\xi}
{\big[(\gm_2-\gm_1)\xi+\gm_1\big]L_1(\xi)+\big[(\gm_2\pi_2-\gm_1\pi_1)\xi+\gm_1\pi_1\big]L_2(\xi)},
\eeq
with $\xi$ being the unknown.
The Newton-Raphson iteration can be used to solve \eqref{eq:f}.

\vspace{2mm}
\subsection{Bound preserving of the volume fraction}\label{subsec:bd-preserv}
It is trivial to check that $f$ defined in \eqref{eq:def-f} is monotonically decreasing with respect to $\tilde\al$ and monotonically increasing with respect to $\xi$, provided that $(\xi,\tilde\al)\in[0,1]\times[0,1-\theta]$. In addition, $f(0;0)=f(1;1-\theta)=0$. So the necessary and sufficient condition for the nonlinear equation $f(\xi;\tilde\al)=0$ possessing a unique solution in the interval $\xi\in[0,1]$ is that $0\leq\tilde{\al}\leq1-\theta$.
In what follows, a restriction on the time step $\Delta t$ is established to ensure the bound-preserving condition $0\leq\tilde{\al}\leq1-\theta$.

By dropping the subscript $i$ and the superscript $n$ indicating  the cell index and the time step, simplify notations in \eqref{eq:def-alpha-star} and rewrite it as
\beq\label{eq:alpha-star-1}
\bga{l}
\tilde{\al}=\bar{\al}-{\Dt}\left(
\beta+\dfr{\Dt}{2}\beta_t
\right)
+\dfr{\Dt}{2}
\bar{\al}
\dfr{\omega-1}{\omega-\bar{\al}}
\eta,
\eda
\eeq
where
\beqs
\beta=
\dfr{(u\al)^{n,*}_{i+\frac 12}-(u\al)^{n,*}_{i-\frac 12}}{\Dx}
, \ \
\beta_t=
\dfr{[(u\al)_t]^{n,*}_{i+\frac 12}-[(u\al)_t]^{n,*}_{i-\frac 12}}{\Dx},
\eeqs
and
\beq\label{eq:def-source-omega}
\omega=\dfr{\rho_1{c_1}^2}{\rho_1{c_1}^2-\rho_2{c_2}^2}.
\eeq
From \eqref{eq:alpha-star-1}, the inequality $\tilde{\al}\geq0$ derives
\beq\label{ineq:lower}
\dfr{\beta_t}{2}\Dt^2
+\left(
\beta
-\dfr{\bar{\al}}{2}
\dfr{\omega-1}{\omega-\bar{\al}}
\eta
\right)
\Dt
\leq\bar{\al}.
\eeq
By basic calculations, the inequality $\tilde{\al}\leq1-\theta$ derives
\beq\label{ineq:upper-2}
\dfr{\hat\beta_t}{2}\Dt^2
+\left(
\hat\beta
-\dfr{\hat\al}{2}\dfr{\hat\omega-1}{\hat\omega-\hat\al}\eta
\right)
\Dt
\leq\hat\al,
\eeq
where $\hat\beta=\eta-\beta$, $\hat\beta_t=\eta_t-\beta_t$, $\hat\al=1-\al$, and $\hat\omega=1-\omega$.
The two inequalities can be treated in the same way since they have the same form. 

Take \eqref{ineq:lower} for example,
whose left hand side is a quadratic function of $\Delta t$. It is easy to find a $\Delta t_0$ such that \eqref{ineq:lower} holds for $\forall\Dt\leq\Dt_0$. 
Therefore, besides the CFL condition, the time step $\Delta t$ is additionally restricted due to inequality \eqref{ineq:lower}.
However, $\Dt_0$ approaches $0$ as $\bar{\al}\rightarrow0^+$ if $\beta-\frac{\bar{\al}}{2}\frac{\omega-1}{\omega-\bar{\al}}\eta>0$. So a cut-off treatment is employed here,
which means that in practical computations, the additional time step restriction according to \eqref{ineq:lower} is only applied in cells where $\bar{\al}>10^{-6}$.
Similarly, the time step is also restricted according to inequality \eqref{ineq:upper-2} in cells where $\bar{\al}<1-10^{-6}$.

\vspace{2mm}
\section{Acoustic Generalized Riemann Problem Solver}\label{sec:grp}
The objective of this section is to develop the GRP solver for the Kapila model under the acoustic assumption, which serves as the building block of the finite volume scheme proposed in the previous section.
There are two versions of the GRP solver: the acoustic version and the nonlinear version. 
The acoustic solver, which linearizes the governing equations, is suitable when the waves emerging from the initial discontinuity are relatively weak \cite{menshov-1991,Toro-ADER}.
Conversely, in the presence of strong waves, the nonlinear GRP solver is required.
The development of the nonlinear GRP solver for the Kapila model is left for a future investigation \cite{grp-kapila}.
To make the present paper self-containing, we develop the GRP solver under the acoustic assumption in this section.

By setting $t^n=0$ and shifting $\bx_{i+\frac 12,j}$ to the origin $\bx_0=(0,0)$, the GRP \eqref{eq:grp-def-local} becomes
\beql{eq:grp-def}
\bga{l}
\dfr{\pt\bU}{\pt t}+\dfr{\pt\mathbf{F}(\bU)}{\pt x}+\dfr{\pt\mathbf{G}(\bU)}{\pt y}=\bK(\bU,\nabla\bU),\\[2mm]
\bU(\bx,t=0)=\left\{
\bga{ll}
\bU_L(\bx), & x<0,\\
\bU_R(\bx), & x>0.
\eda
\right.
\eda
\eeq
In finite volume schemes, the piecewise smooth initial data are obtained by the space data reconstruction procedure in cells that are adjacent to the initial discontinuity $\{\bx:x=0\}$.

\vspace{2mm}
\subsection{Solve the associated Riemann problem}\label{subsec:ass-rp}
Following Ref.~\onlinecite{grp-06}, the first step of solving the GRP  \eqref{eq:grp-def} is to solve the associated Riemann problem
\beql{eq:ass-rp-def}
\bga{l}
\dfr{\pt\bU^\text{ass}}{\pt t}+\dfr{\pt\bF(\bU^\text{ass})}{\pt x}
=\bK^\text{ass}(\bU^\text{ass},\nabla\bU^\text{ass}),\\
\bU^\text{ass}(\bx,t=0)=\left\{
\bga{ll}
\bU_L^\text{ass}=\bU_L(\bx_0), & x<0,\\
\bU_R^\text{ass}=\bU_R(\bx_0), & x>0,
\eda
\right.
\eda
\eeq
where
\beqs
\mathbf{K}^\text{ass}=[0,0,0,0,0,\al_1K\dfr{\pt u}{\pt x}]^\top.
\eeqs

A concern that arises in this situation is the ambiguity of the Hugoniot relation determining the shock behavior, which is a result of the governing equations being genuinely non-conservative.
Here we use the one proposed in Ref.~\onlinecite{saurel-shock-jump},
\beql{eq:hugoniot}
\bga{l}
\zeta_k=\zeta_k^0, \\
\rho(u-\sigma)=\rho^0(u^0-\sigma)\triangleq m, \\
p-p^0=m^2(\tau-\tau^0), \\
e_k-e_k^0+\dfr{p+p^0}{2}(\tau_k-\tau_k^0)=0,
\eda
\eeq
where $\sigma$ is the shock speed, $\tau=\frac{1}{\rho}$ is the specific volume of the mixture, and $q^0$ denotes the pre-shock value of a certain physical quantity $q$. The advantages of the Hugoniot relation \eqref{eq:hugoniot} are highlighted as follows \cite{saurel-shock-jump,saurel-projection-II}: (i) The conservation of the total energy is preserved. (ii) It resembles to single-phase shock relations. (iii) It preserves the positivity of volume fractions. (iv) It is symmetric with respect to all phases. (v) It is tangent to the mixture isentrope, which is important for the well-posedness of the Riemann problem. (vi) It agrees with a variety of experimental measurements.

Based on Hugoniot relation \eqref{eq:hugoniot}, the IVP \eqref{eq:ass-rp-def} can be solved analytically by the exact Riemann solver \cite{saurel-projection-II}. The Riemann solution is therefore obtained and denoted as
\begin{equation}\label{eq:rs}
\bU^*=R^A(0,\bU_L^\text{ass},\bU_R^\text{ass})
=\displaystyle\lim_{t\rightarrow0^+}\bU^\text{ass}(\bx_0,t).
\end{equation}

\vspace{2mm}
\begin{rem}
Note that the initial data used in \eqref{eq:ass-rp-def} are a piecewise constant functions by evaluating $\bU_{L/R}$ at $\bx_0$, which means that the spatial variation of the flow field is neglected in numerical computations if the scheme is developed solely based on the Riemann solution. 
Moreover, by considering $\mathbf{K}^\text{ass}$ instead of $\bK$ as the source term, only the normal component (the $x$-derivative) of the velocity divergence is taken into account and the transversal component is dropped. Therefore, the compaction or expansion effect of the flow field is not fully included in the Riemann solution. 
\end{rem}

\vspace{2mm}
\begin{rem}
The ambiguity of the Hugoniot relation also arises at the discrete level, since numerical solutions of non-conservative hyperbolic PDEs depend on the way to discretize governing equations.
Nevertheless, the numerical experiments conducted in Section \ref{sec:numer} indicate that the treatment employed by the finite volume scheme \eqref{eq:scheme} produces numerical results that agree with exact solutions determined by the Hugoniot relation \eqref{eq:hugoniot}.
\end{rem}

\vspace{2mm}
\subsection{Solve the acoustic generalized Riemann problem}\label{subsec:grp}

Regarding primitive variables as unknowns, rewrite the governing equations in the quasi-linear form
\begin{subequations}\label{eq:govern-quasi}
\begin{align}
&\dfr{\pt\bV}{\pt t}+\bA\dfr{\pt\bV}{\pt x}+\bB\dfr{\pt\bV}{\pt y}=0,    \label{eq:govern-quasi-1}\\[3mm]
&\dfr{\pt\al_1}{\pt t}+u\dfr{\pt\al_1}{\pt x}+v\dfr{\pt\al_1}{\pt y}=(K-\alpha_1) \ 
\left(\dfr{\pt u}{\pt x}+\dfr{\pt v}{\pt y}\right),    \label{eq:govern-quasi-2}
\end{align}
\end{subequations}
where
\beql{eq:govern-quasi-def}
\bga{l}
\bV=\left[
\bga{c}
\zeta_1\\
\rho\\
u\\
v\\
p
\eda
\right],\ \ 
\bA=\left[
\bga{ccccc}
u & 0 & 0 & 0 & 0\\
0 & u & \rho & 0 & 0\\
0 & 0 & u & 0 & \dfr{1}{\rho}\\
0 & 0 & 0 & u & 0\\
0 & 0 & \rho c^2 & 0 & u
\eda
\right], \ \
\bB=\left[
\bga{ccccc}
v & 0 & 0 & 0 & 0\\
0 & v & 0 & \rho & 0\\
0 & 0 & v & 0 & 0\\
0 & 0 & 0 & v & \dfr{1}{\rho}\\
0 & 0 & 0 & \rho c^2 & v
\eda
\right].
\eda
\eeq
By regarding \eqref{eq:govern-quasi} as governing equations, the IVP \eqref{eq:grp-def} becomes
\beql{eq:grp-quasi-def}
\bga{l}
\dfr{\pt\bV}{\pt t}+\bA\dfr{\pt\bV}{\pt x}+\bB\dfr{\pt\bV}{\pt y}=0, \\[2.5mm]
\dfr{\pt\al_1}{\pt t}+u\dfr{\pt\al_1}{\pt x}+v\dfr{\pt\al_1}{\pt y}=(K-\alpha_1) \ 
\left(\dfr{\pt u}{\pt x}+\dfr{\pt v}{\pt y}\right),\\[2.5mm]
\bV(\bx,t=0)=\left\{
\bga{ll}
\bV_L(\bx), & x<0,\\
\bV_R(\bx), & x>0,
\eda
\right.\\[2mm]
\al_1(\bx,t=0)=\left\{
\bga{ll}
(\al_1)_L(\bx), & x<0,\\
(\al_1)_R(\bx), & x>0.
\eda
\right.
\eda
\eeq
The initial values of $\bV$ and $\al_1$ are obtained by transforming the reconstructed data of $\bU$ into those of primitive variables.
By noting that \eqref{eq:govern-quasi-1} is independent of \eqref{eq:govern-quasi-2}, the acoustic GRP solver \cite{grp-06} can be directly applied to solve the IVP of \eqref{eq:govern-quasi-1}.

With the Riemann solution $\bU^*$ obtained in \eqref{eq:rs}, fix the coefficient matrices as $\bA^*=\bA(\bU^*)$ and $\bB^*=\bB(\bU^*)$. By the characteristic decomposition,
\beql{eq:grp-res}
\bga{l}
\left(\dfr{\pt\bV}{\pt t}\right)^*=
-\bR^*\mathbf\Lambda_+^{*}\bL^{*}\left.\dfr{\pt\bV_L}{\pt x}\right|_{\bx=\mathbf{x}_0}
-\bR^*\mathbf\Lambda_-^{*}\bL^{*}\left.\dfr{\pt\bV_R}{\pt x}\right|_{\bx=\mathbf{x}_0}\\[4mm]
\qquad\qquad\quad
-\bR^*\bI_+^{*}\bL^{*}\bB^*\left.\dfr{\pt\bV_L}{\pt y}\right|_{\bx=\mathbf{x}_0}
-\bR^*\bI_-^{*}\bL^{*}\bB^*\left.\dfr{\pt\bV_R}{\pt y}\right|_{\bx=\mathbf{x}_0}.
\eda
\eeq
The matrices $\bR^*$, $\bL^*$ and $\mathbf\Lambda^*$ are given by the decomposition of $\bA^*$ as
\beqs
\bA^*=\bR^*\mathbf\Lambda^*\bL^{*},
\eeqs
where $\bL^*$ is the inverse of $\bR^*$.
The matrix $\bR=[\bR_1,\bR_2,\bR_3,\bR_4,\bR_5]$ is composed of the right eigenvectors of $\bA$, which are
\begin{equation*}
\bga{c}
\bR_1=[0, \dfr{1}{c^2}, -\dfr{1}{\rho c}, 0, 1]^\top, \ \
\bR_2=[0, 0, 0, 1, 0]^\top, \ \
\bR_3=[0, 1, 0, 0, 0]^\top, \\
\bR_4=[1, 0, 0, v, \dfr{1}{\rho}]^\top, \ \
\bR_5=[0, \dfr{1}{c^2}, \dfr{1}{\rho c}, 0, 1]^\top.
\eda
\end{equation*}
Corresponding eigenvalues of $\bA$ are $\lambda_1=u-c$, $\lambda_2=\lambda_3=\lambda_4=u$, and $\lambda_5=u+c$. Furthermore,
\beqs
\bga{l}
\mathbf\Lambda_+=\text{diag}\{\max(\lambda_i,0)\}, \ \
\mathbf\Lambda_-=\text{diag}\{\min(\lambda_i,0)\},\\
\bI_+=\text{diag}\{\max(\text{sign}(\lambda_i),0)\}, \ \
\bI_-=\text{diag}\{\max(\text{sign}(-\lambda_i),0)\}.
\eda
\eeqs

The last step is to compute the time derivative of the volume fraction $\al_1$. Instead of directly solving the IVP of \eqref{eq:govern-quasi-2}, in what follows, we derive the value of $(\pt\al_1/\pt t)^*$ from $(\pt\bV/\pt t)^*$ obtained in \eqref{eq:grp-res}.
The time derivative of the phase density is
\beqs
\bga{l}
\dfr{\pt\rho_1}{\pt t}
=\left[\left(\dfr{\pt p}{\pt\rho_1}\right)_{s_1}\right]^{-1}
\left[
\dfr{\pt p}{\pt t}
-\left(\dfr{\pt p}{\pt s_1}\right)_{\rho_1}
\dfr{\pt s_1}{\pt t}
\right]\\[3mm]
\qquad
=\dfr{1}{{c_{1}}^2}\left[
\dfr{\pt p}{\pt t}
+u\left(\dfr{\pt p}{\pt s_1}\right)_{\rho_1}
\dfr{\pt s_1}{\pt x}
\right].
\eda
\eeqs
Multiplying the above equation by $\tau_1$ and taking the limit $t\rightarrow0^+$ lead to
\beq\label{eq:grp-res-taurhot}
\bga{l}
\left(\tau_1\dfr{\pt\rho_1}{\pt t}\right)^*
=\dfr{1}{\rho_{1}^*{c_{1}^*}^2}\left[
\left(\dfr{\pt p}{\pt t}\right)^*
+u^*\left(\dfr{\pt p}{\pt s_1}\right)_{\rho_1}^*
\left(\dfr{\pt s_1}{\pt x}\right)^*
\right]
,
\eda
\eeq
where $(\frac{\pt p}{\pt t})^*$ is obtained in \eqref{eq:grp-res}. For stiffened gases, $\rho_{1}^*{c_{1}^*}^2=\gm_1(p^*+\pi_1)$.
By the fact that $\al_1\tau=\zeta_1\tau_1$,
\beqs
\bga{l}
\dfr{1}{\tau}d\tau+\dfr{1}{\al_k}d\al_k
=\dfr{1}{\tau_k}d\tau_k+\dfr{1}{\zeta_k}d\zeta_k\\
\qquad\qquad\qquad \
=-\tau_kd\rho_k+\dfr{1}{\zeta_k}d\zeta_k.
\eda
\eeqs
Therefore, by approaching the limit $t\rightarrow 0^+$,
\beq\label{eq:grp-res-vof}
\bga{l}
\left(\dfr{\pt\al_1}{\pt t}\right)^*
=\al_1^*\left(
-\tau_1\dfr{\pt \rho_1}{\pt t}+\dfr{1}{\zeta_1}\dfr{\pt \zeta_1}{\pt t}+\tau\dfr{\pt\rho}{\pt t}
\right)^*.
\eda
\eeq
The acoustic GRP solver is finalized.

\vspace{2mm}
\begin{rem}
It is illustrated by Lei and Li \cite{Lei-Li-2018-AMM} that the GRP solver improves corresponding numerical schemes to effectively resolve multi-dimensional effects of the flow field. 
This advantage is of particular importance for the current scheme developed for the Kapila model, since the compaction or expansion effect, which is indeed a multi-dimensional effect, plays a key role in the volume fraction equation.
By substituting governing equations \eqref{eq:govern} into \eqref{eq:grp-res-vof}, the time derivative of $\alpha_1$ is expressed as
\beql{eq:grp-res-al}
\bga{l}
\displaystyle
\left(\dfr{\pt\al_1}{\pt t}\right)^*
=\lim_{t\rightarrow0^+}\left.\dfr{\pt\al_1}{\pt t}\right|_{(\bx=\mathbf{0},t)}\\[3.5mm]
\qquad\qquad
=\displaystyle\lim_{t\rightarrow0^+}
\left.\left[-u\dfr{\pt\al_1}{\pt x}-v\dfr{\pt\al_1}{\pt y}+(K-\alpha_1) \ 
\left(\dfr{\pt u}{\pt x}+\dfr{\pt v}{\pt y}\right)\right]\right|_{(\bx=\mathbf{0},t)},
\eda
\eeq
which is the Lax-Wendroff procedure applied to the volume fraction equation \eqref{eq:vof}.
Compared with the Riemann solution obtained by solving the associated Riemann problem \eqref{eq:ass-rp-def}, the divergence of the velocity field, which reflects the compaction or expansion effect of the flow field, is fully incorporated into the time derivative \eqref{eq:grp-res-vof}.
\end{rem}

\vspace{2mm}
\section{Numerical Experiments}\label{sec:numer}
This section conducts both one- and two-dimensional numerical experiments.
The CFL number is taken to be $CFL=0.6$ except for two-dimensional cases
where $CFL=0.45$ for Example 5 and $CFL=0.25$ for Example 6. 
The nonlinear GRP solver is used for the water-air shock tube problem in Example 2 since strong nonlinear waves emanating from the initial discontinuity activate severe interactions between phases.
The acoustic GRP solver is used for other cases.
In all cases, the parameter for the linear data reconstruction \eqref{eq:recon-2d} is $\kappa=1.5$.

\vspace{2mm}
\n{\bf Example 1. Accuracy test.}
We first perform the accuracy test by considering the evolution of the mixture of water and air. 
Both water and air follow the stiffened gas EOS \eqref{eq:eos-sg}. The parameters for water are taken to be $\gm_1=4.4$ and $\pi_1=6000$, while those for air are taken to be $\gm_2=1.4$ and $\pi_2=0$.
The initial phase density of water is defined as a periodic function
\beqs
\rho_1(x,0)=20+2\sin(2\pi x).
\eeqs
The pressure of the mixture is defined according to the isentropic condition of water as
\beqs
p(x,0) = S_1 {\rho_1}^{\gm_1} - \pi_1,
\eeqs
where $S_1=0.05$.
The air is also initially isentropic, making its phase density
\beqs
\rho_2(x,0) = \left(\dfr{p(x,0)+\pi_2}{S_2}\right)^{\frac{1}{\gm_2}},
\eeqs
where $S_2=5000$.
The volume fraction of water is determined by the relation
\beqs
\bga{l}
\al_1(x,0)\rho_1(x,0) = \zeta(x,0)\rho(x,0) \\
\qquad\qquad\qquad\ \ \
= \zeta(x,0)\big[\al_1(x,0)\rho_1(x,0) + (1-\al_1(x,0))\rho_2(x,0)\big].
\eda
\eeqs
The mass fraction of water is taken to be constant as $\zeta_1(x,0)=0.992$.
The mixture entropy is defined as
\beq\label{eq:total-entropy}
S=\zeta_1h_1(S_1)+\zeta_2h_2(S_2),
\eeq
where $h_k(S_k)={\zeta_k}^{\gamma_k}S_k$ is the entropy of each phase $k$ for $k=1,2$. The computational domain is set to be $0\leq x\leq 1$, and the periodic boundary condition is used at both ends.

The exact solution for $\rho$, $u$, or $p$ cannot be analytically expressed. However, the solution remains isentropic as long as the flow is smooth.
So the errors with respect to the water entropy $h_1(S_1)$ and the air entropy $h_2(S_2)$ are defined as
\begin{equation*}
\mathcal{E}_k(x,t)=h_k\left(\dfr{p(x,t)}{[\rho_k(x,t)]^{\gamma_k}}\right)-h_k(S_k), \ \ k=1,2,
\end{equation*}
and the error with respect to the mixture entropy is defined as
\begin{equation*}
\mathcal{E}(x,t)=\sum_{k=1,2}\zeta_k(x,t)h_k\left(\dfr{p(x,t)}{[\rho_k(x,t)]^{\gamma_k}}\right) - S.
\end{equation*}
The errors $\mathcal{E}_1$ and $\mathcal{E}_2$ are evaluated at the computational time $T=5\times 10^{-3}$. Their $L_1$ and $L_\infty$ norms and corresponding convergence rates are listed in Tab. \ref{tab:accuracy-test-1d-fractional}. Those of the mixture entropy $S$ are listed in Tab. \ref{tab:accuracy-test-1d-total}. It is observed that the proposed scheme attains the designed second-order accuracy.

\begin{table*}[!htbp]
  \centering
  \caption[small]{$L_1$ and $L_\infty$ errors of the water entropy $h_1(S_1)$ and the air entropy $h_2(S_2)$, and corresponding convergence rates.}
\footnotesize
  \begin{tabular}{|c|cccc|cccc|}
    \hline
                     &    \multicolumn{4}{|c|}{error of water entropy}     &  \multicolumn{4}{|c|}{error of  air entropy}   \\\cline{2-9}
   mesh size    & $L_1$ error &  order   & $L_\infty$ error & order  & $L_1$ error &  order   & $L_\infty$ error & order  \\\hline
 1/20 & 3.59e-04 &  -   & 8.55e-04 &  -   & 7.95e-04 &  -   & 1.29e-03 &  -   \\
 1/40 & 9.48e-05 & 1.92 & 2.41e-04 & 1.83 & 1.94e-04 & 2.04 & 3.20e-04 & 2.01 \\
 1/80 & 2.45e-05 & 1.95 & 6.48e-05 & 1.89 & 4.78e-05 & 2.02 & 7.67e-05 & 2.06  \\
 1/160 & 6.14e-06 & 2.00 & 1.66e-05 & 1.96 & 1.21e-05 & 1.98 & 1.90e-05 & 2.02 \\
 1/320 & 1.50e-06 & 2.03 & 3.81e-06 & 2.13 & 3.48e-06 & 1.80 & 6.49e-06 & 1.55  \\
 1/640 & 3.65e-07 & 2.04 & 9.58e-07 & 1.99 & 9.00e-07 & 1.95 & 1.74e-06 & 1.90 \\\hline 
  \end{tabular}
  \label{tab:accuracy-test-1d-fractional}
\end{table*}

\begin{table*}[!htbp]
  \centering
  \caption[small]{$L_1$ and $L_\infty$ errors of the mixture entropy $S$, and corresponding  convergence rates.}
\footnotesize
  \begin{tabular}{|c|cccc|}
    \hline
 \multicolumn{5}{|c|}{error of mixtue entropy}  \\\hline
   mesh size    & $L_1$ error &  order   & $L_\infty$ error & order   \\\hline 
 1/20 & 5.50e-04 &  -   & 9.93e-04 &  -    \\
 1/40 & 1.36e-04 & 2.02 & 2.74e-04 & 1.86  \\
 1/80 & 3.38e-05 & 2.00 & 6.89e-05 & 1.99  \\
 1/160 & 8.56e-06 & 1.98 & 1.73e-05 & 1.99  \\
 1/320 & 2.27e-06 & 1.91 & 4.84e-06 & 1.84  \\
 1/640 & 5.79e-07 & 1.97 & 1.21e-06 & 2.00  \\\hline 
  \end{tabular}
  \label{tab:accuracy-test-1d-total}
\end{table*}

\vspace{2mm}
\n{\bf Example 2. Water-air shock tube problem.} This case is widely used to test the capability of a numerical scheme to handle the interface separating fluids exhibiting large density and pressure ratios \cite{saurel-projection-II,saurel-jfm,murrone}. Here we consider the density ratio to be $1000:1$.
Let water and air occupy the interval $[0,1]$. The EOSs of water and air are both given by \eqref{eq:eos-sg}. The parameters of the water EOS are taken to be $\gm_1=4.4$ and $\pi_1=6\times 10^8$, while those of the air EOS are taken to be $\gm_2=1.4$ and $\pi_2=0$. 
Initially, the fluid remains at rest with the following initial condition,
\beqs
(\al_1, \ \rho_1, \ \rho_2, \ p)=
\left\{
\bga{lr}
(1, \ 1000, \ 1, \ 10^9), & x<0.7,\\[2mm]
(0, \ 1000, \ 1, \ 10^5), & x>0.7.
\eda
\right.
\eeqs

\begin{figure}[!htb]
\centering
\subfigure[mixture density]{
\includegraphics[width=.46\linewidth]{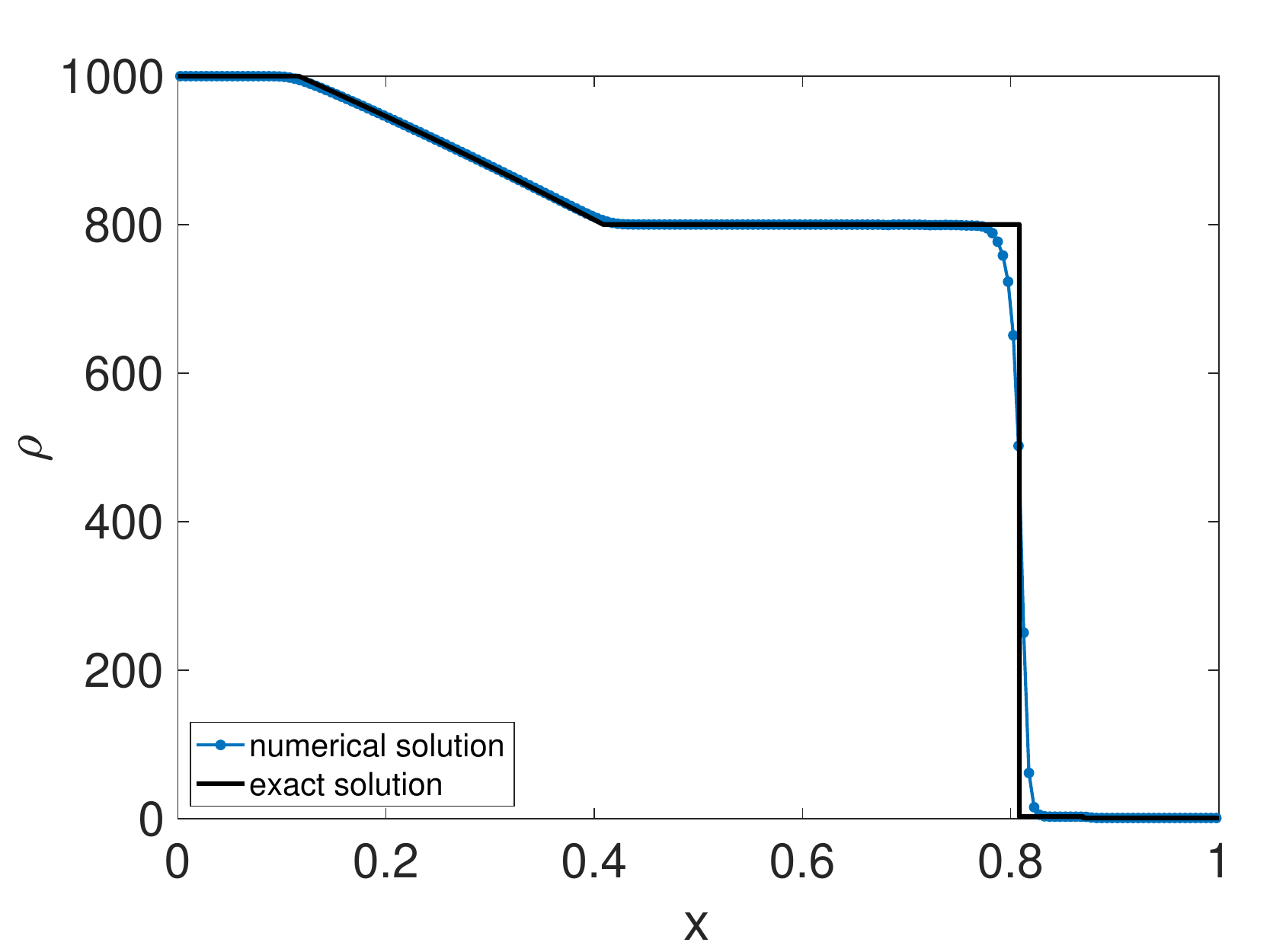}
\label{fig:water-air-shocktube-1000-density}
}
\subfigure[pressure]{
\includegraphics[width=.46\linewidth]{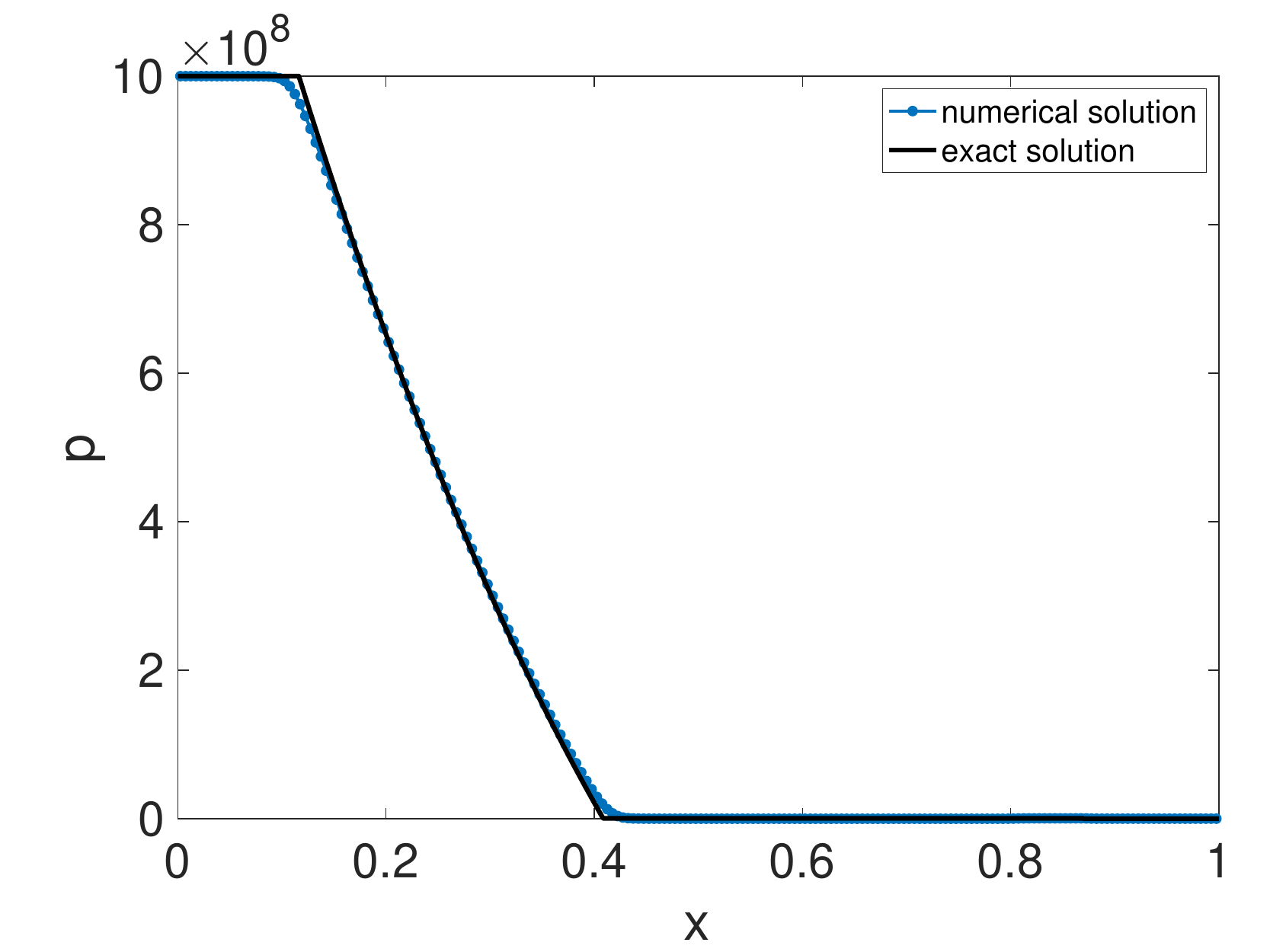}
\label{fig:water-air-shocktube-1000-pressure}
}
\subfigure[gas volume fraction]{
\includegraphics[width=.46\linewidth]{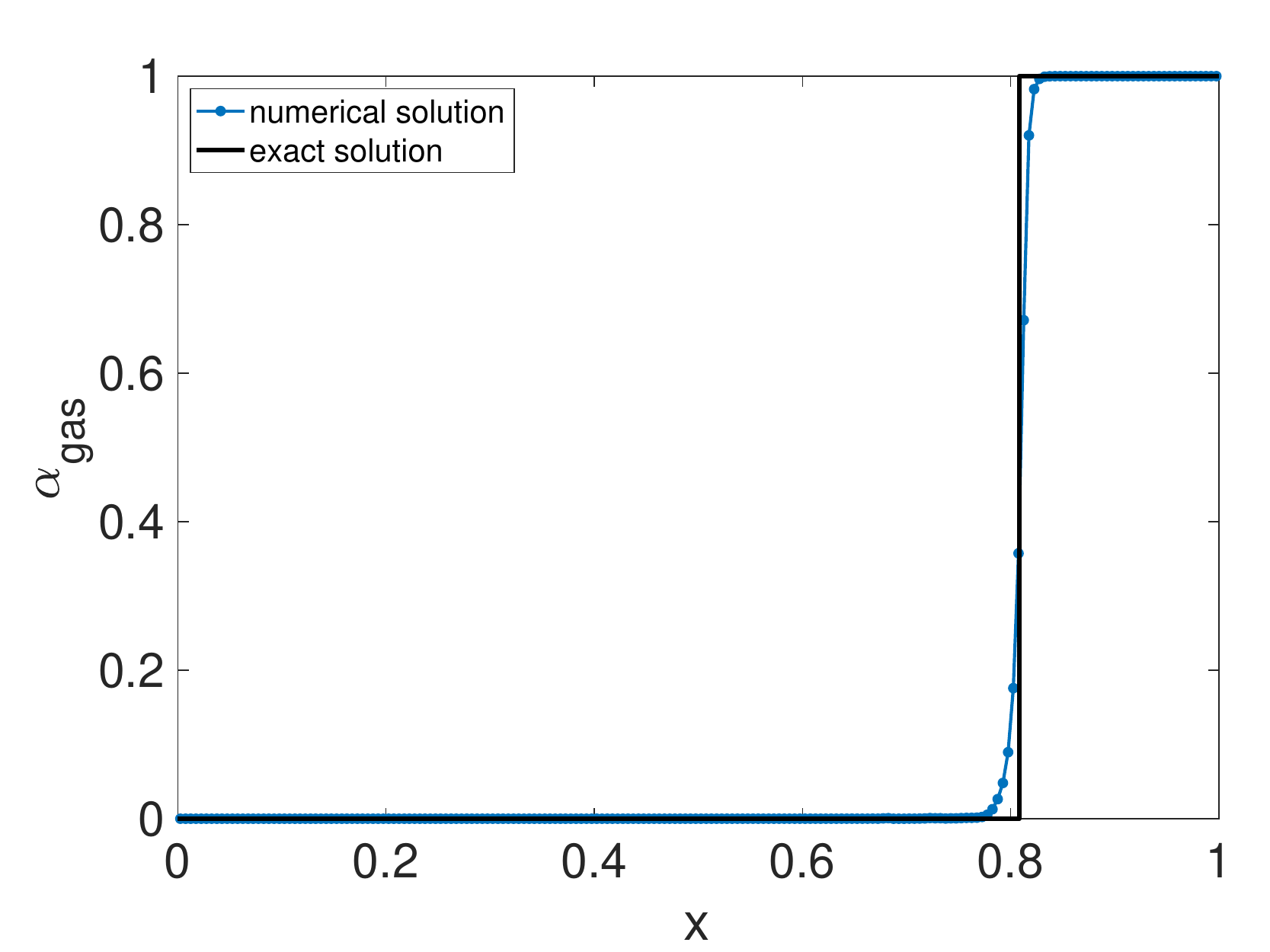}
\label{fig:water-air-shocktube-1000-vof}
}
\subfigure[velocity]{
\includegraphics[width=.46\linewidth]{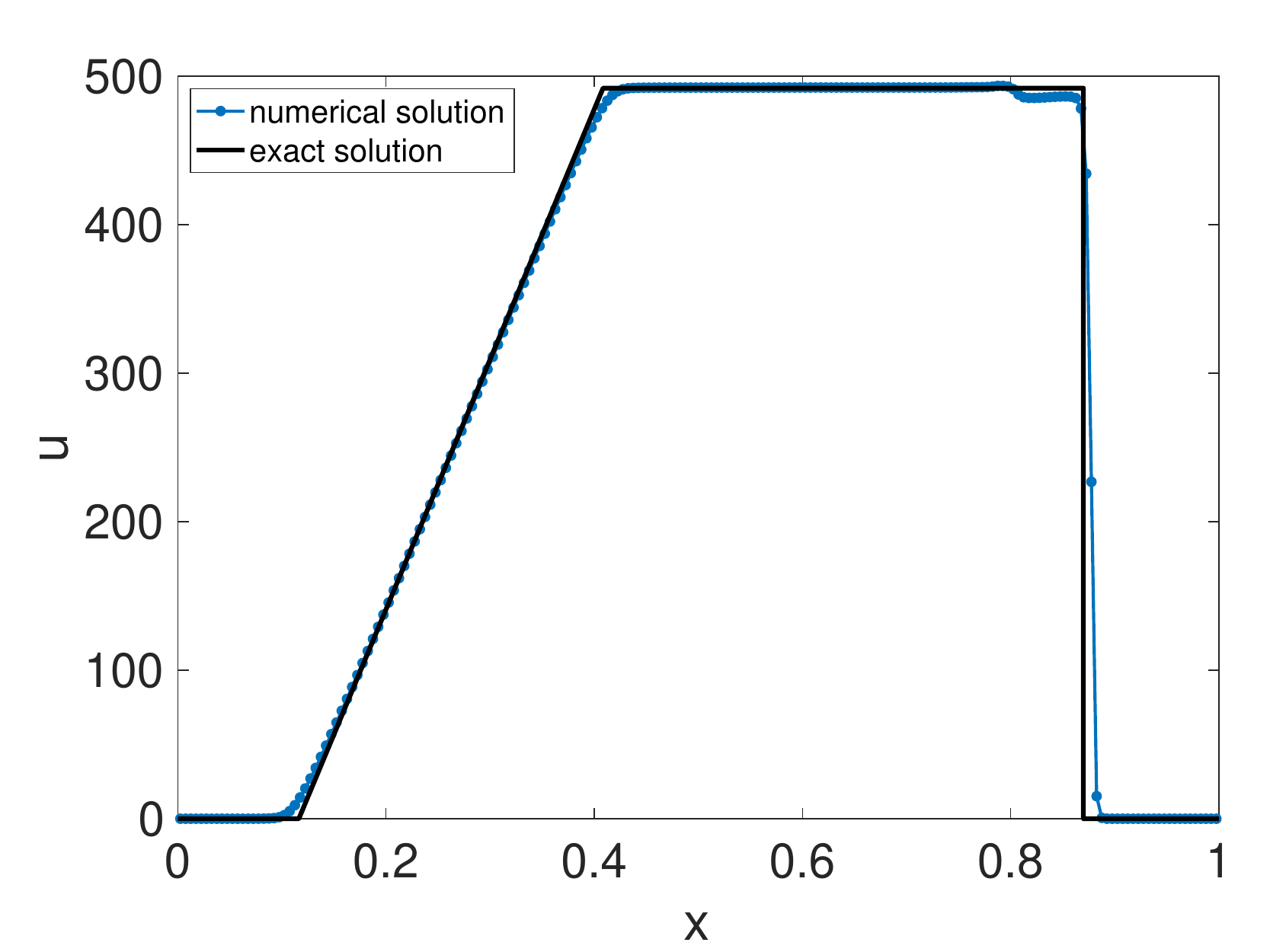}
\label{fig:water-air-shocktube-1000-velocity}
}
\caption{Numerical solutions (dots), in comparison with exact ones (black lines), of the water-air shock tube problem. The mixture density, pressure, gas volume fraction, and velocity are displayed in subfigures.}
\label{fig:water-air-shocktube-1000}
\end{figure}

We perform the simulation with $200$ equally distributed grids to the time $T=2.2\times 10^{-4}$. The numerical results obtained by the proposed finite volume scheme are displayed in Fig. \ref{fig:water-air-shocktube-1000}, which match well with the exact solution.

In this particular test case, it is important to emphasize that the strong rarefaction wave introduces severe thermodynamical variations in space and time. 
Therefore the nonlinear GRP solver becomes indispensable for the numerical scheme to correctly capture the evolution of the flow field.

\vspace{2mm}
\n{\bf Example 3. Two-phase water-air problem.} This is the case considered by Murrone and Guillard \cite{murrone}. Let the mixture of water and air occupy the interval $[0,1]$. The EOSs of water and air are the same as those used in the previous case. The mixture remains at rest initially. The initial condition is given as follows:
\beqs
(\al_1, \ \rho_1, \ \rho_2, \ p)=
\left\{
\bga{lr}
(0.5, \ 1000, \ 50, \ 10^9), & x<0.5,\\[2mm]
(0.5, \ 1000, \ 50, \ 10^5), & x>0.5.
\eda
\right.
\eeqs

\begin{figure}[!htb]
\centering
\subfigure[mixture density]{
\includegraphics[width=.46\linewidth]{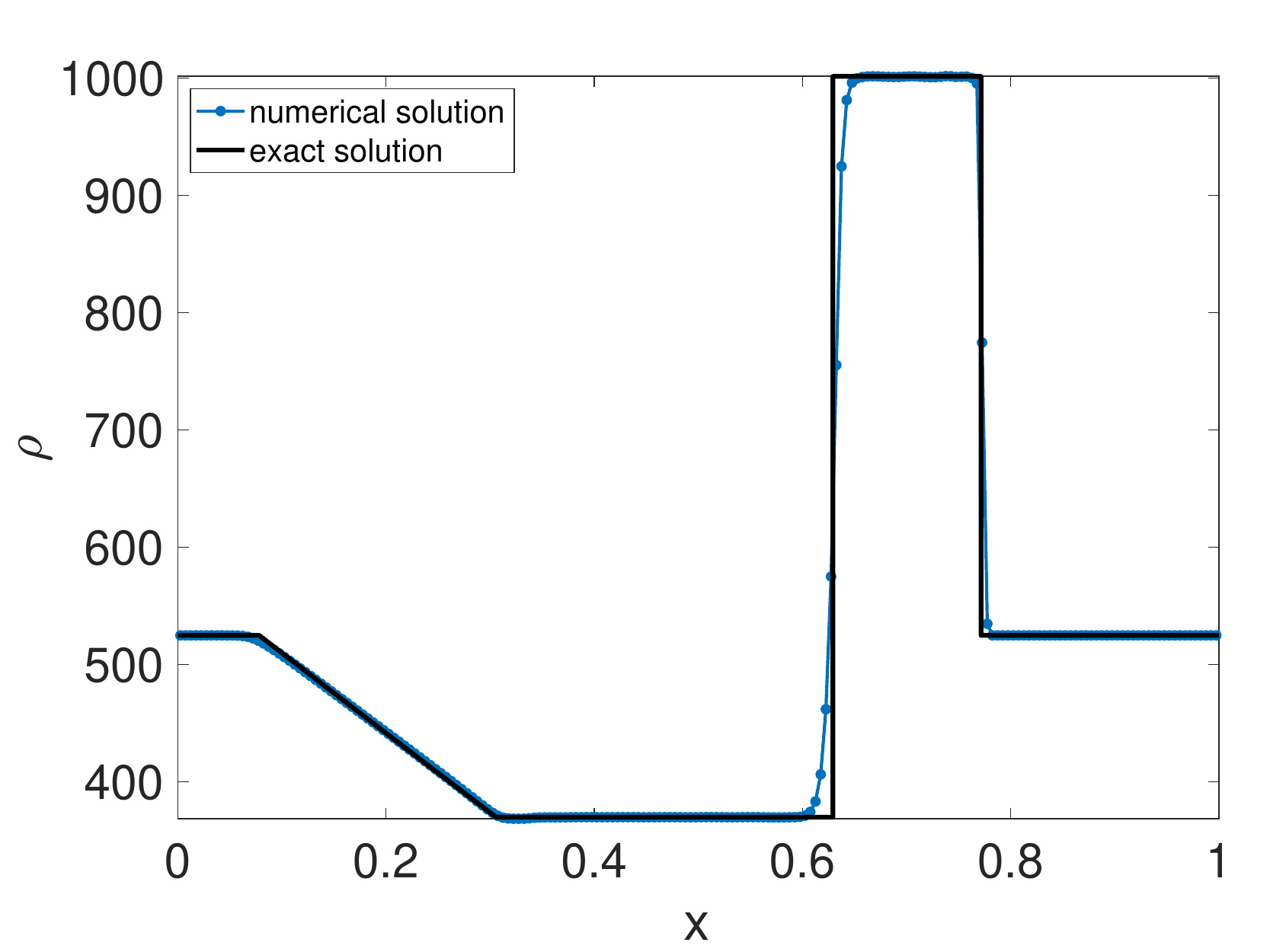}
\label{fig:two-phase-density}
}
\subfigure[pressure]{
\includegraphics[width=.46\linewidth]{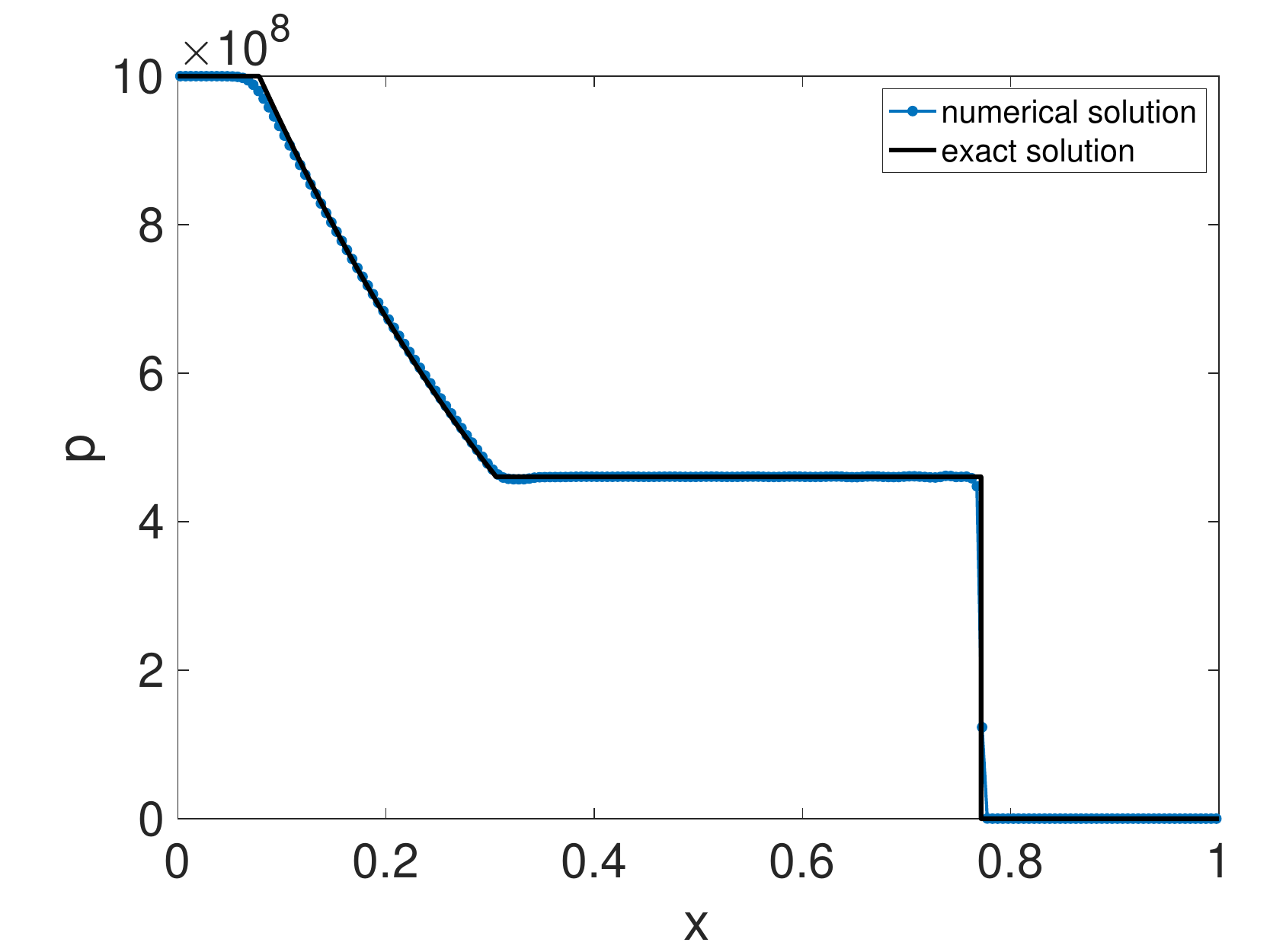}
\label{fig:two-phase-pressure}
}
\subfigure[gas volume fraction]{
\includegraphics[width=.46\linewidth]{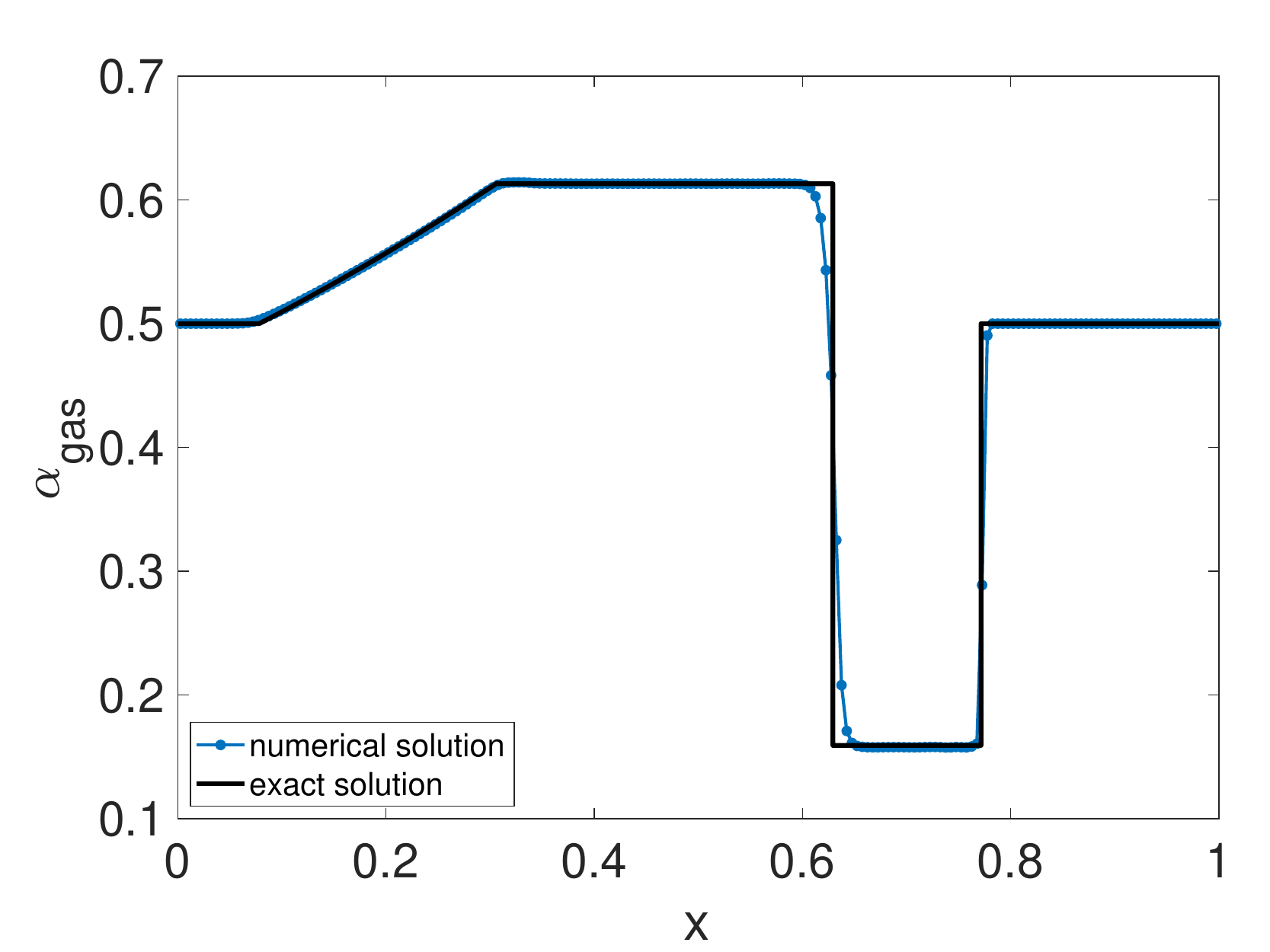}
\label{fig:two-phase-vof}
}
\subfigure[velocity]{
\includegraphics[width=.46\linewidth]{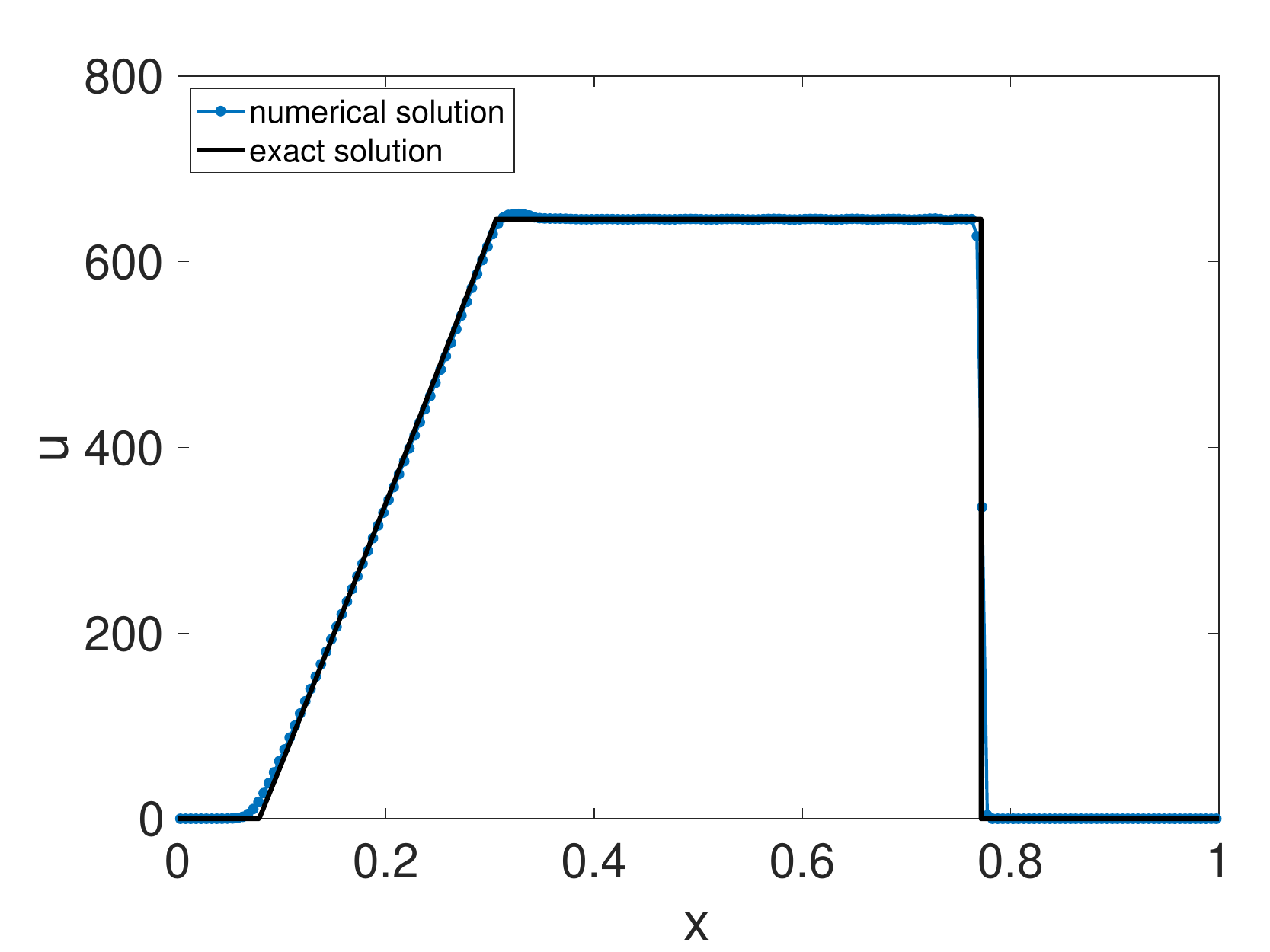}
\label{fig:two-phase-velocity}
}
\caption{Numerical solutions (dots), in comparison with exact ones (black lines), of the two-phase water-air problem. The mixture density, pressure, gas volume fraction, and velocity are displayed in subfigures.}
\label{fig:two-phase}
\end{figure}

Perform the simulation with $200$ equally distributed grids to the time $T=2\times 10^{-4}$. The numerical results obtained by the propose finite volume scheme are displayed in Fig. \ref{fig:two-phase}, and an excellent agreement with exact solutions is observed. In addition, the numerical results match well with those obtained by the seven-equation model presented in Ref.~\onlinecite{murrone}.

\begin{figure}[!htb]
\centering
\subfigure[mixture density]{
\includegraphics[width=.46\linewidth]{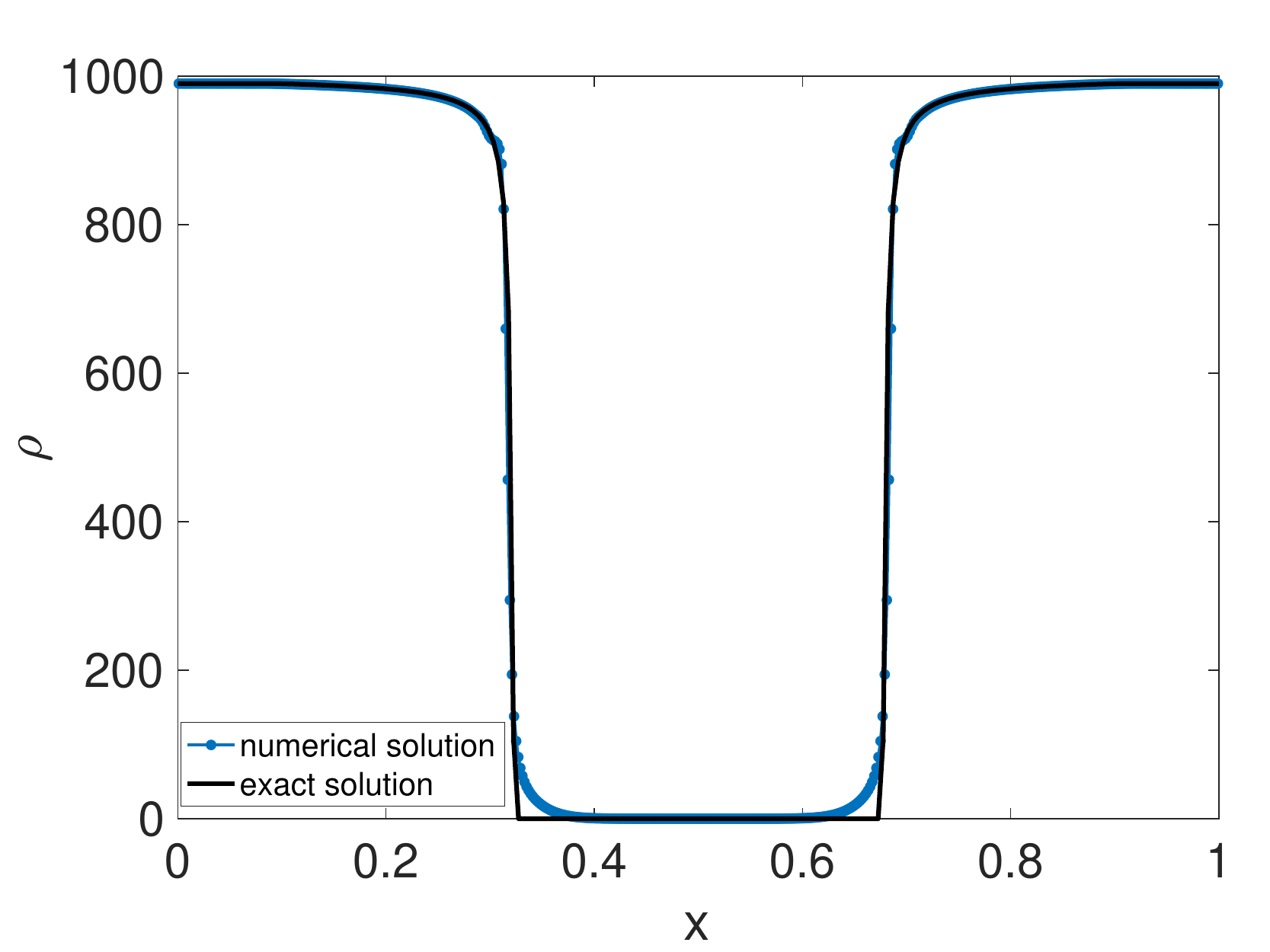}
\label{fig:cavitation-density}
}
\subfigure[pressure]{
\includegraphics[width=.46\linewidth]{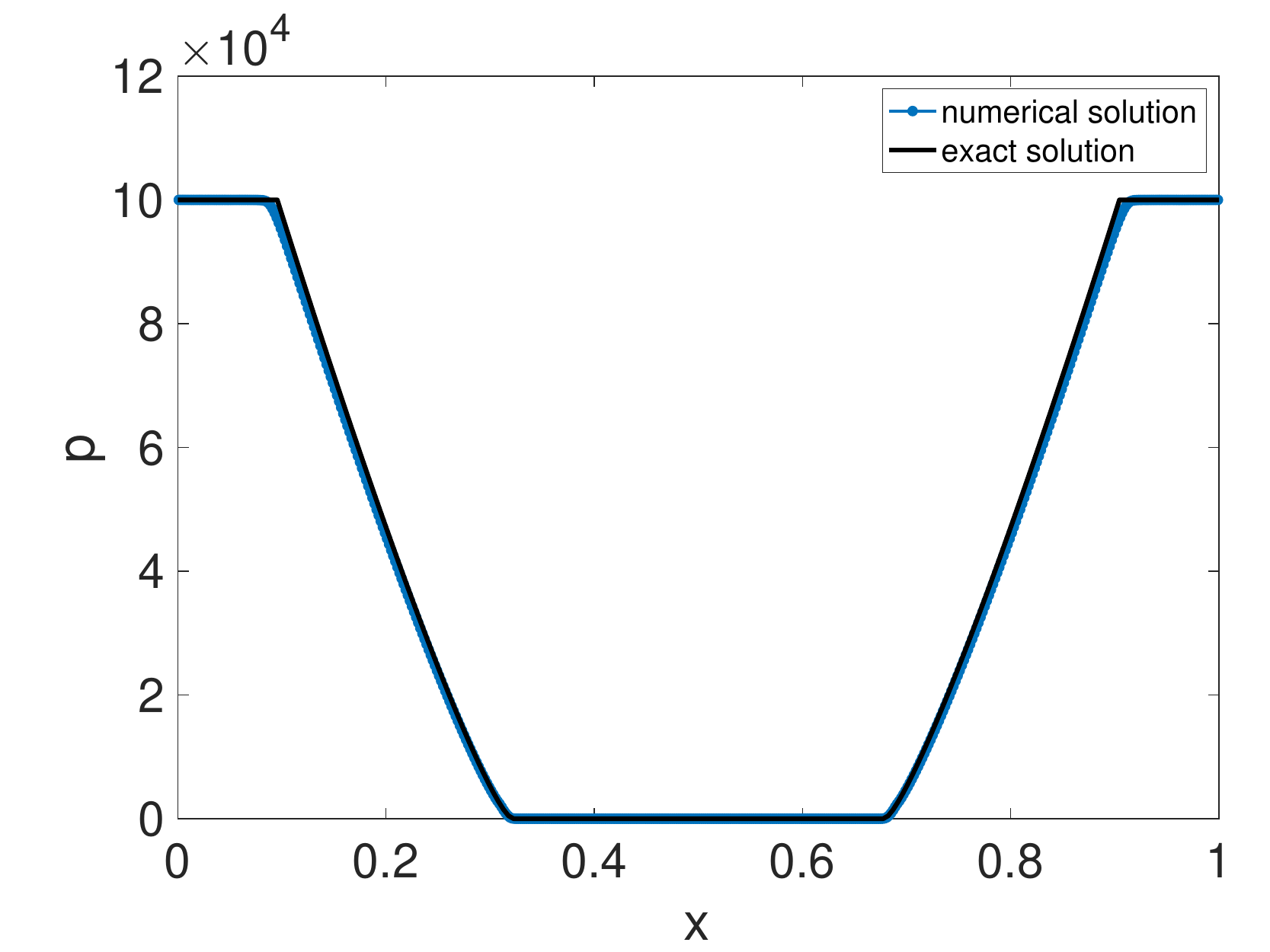}
\label{fig:cavitation-pressure}
}
\subfigure[gas volume fraction]{
\includegraphics[width=.46\linewidth]{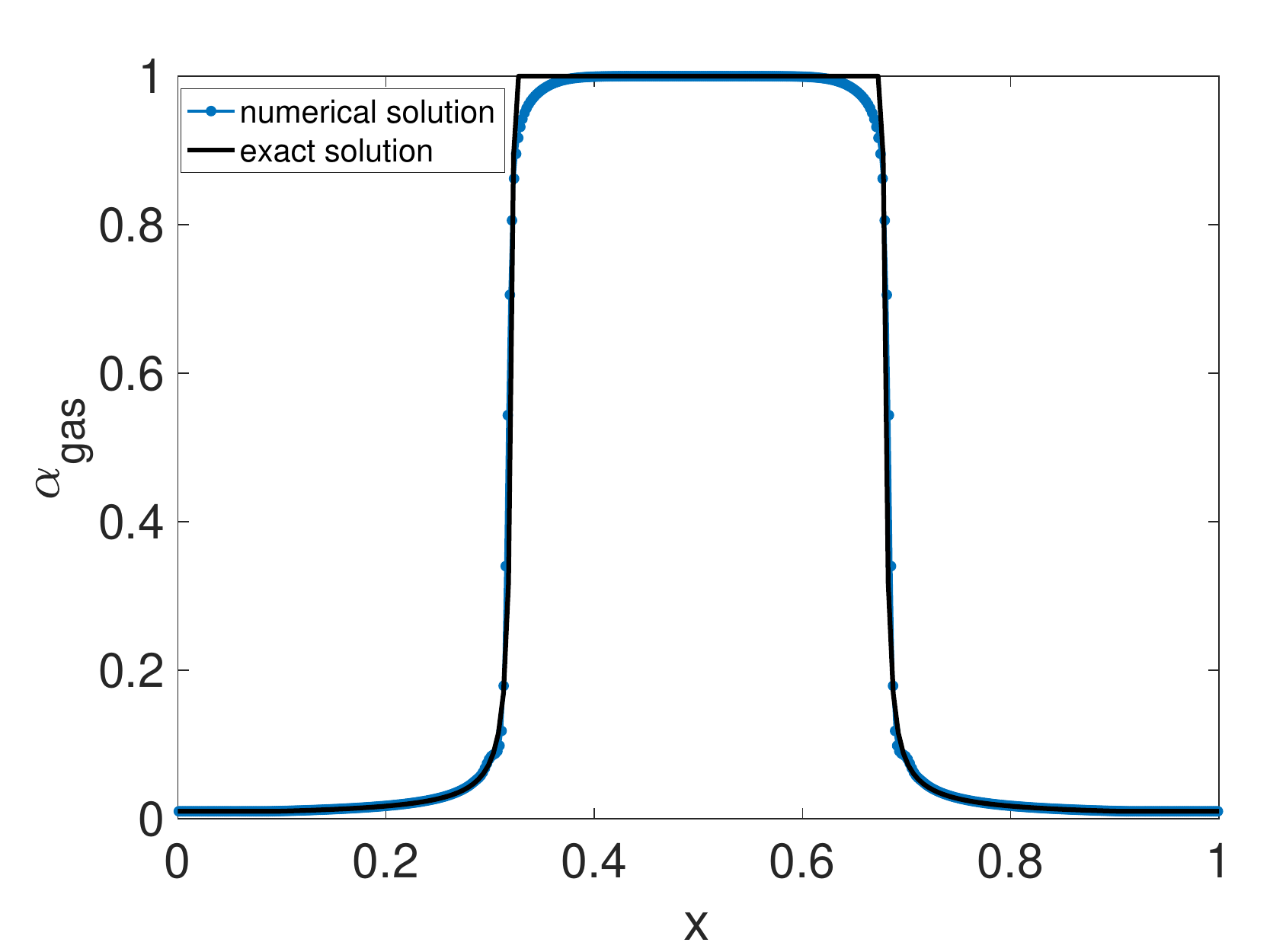}
\label{fig:cavitation-vof}
}
\subfigure[velocity]{
\includegraphics[width=.46\linewidth]{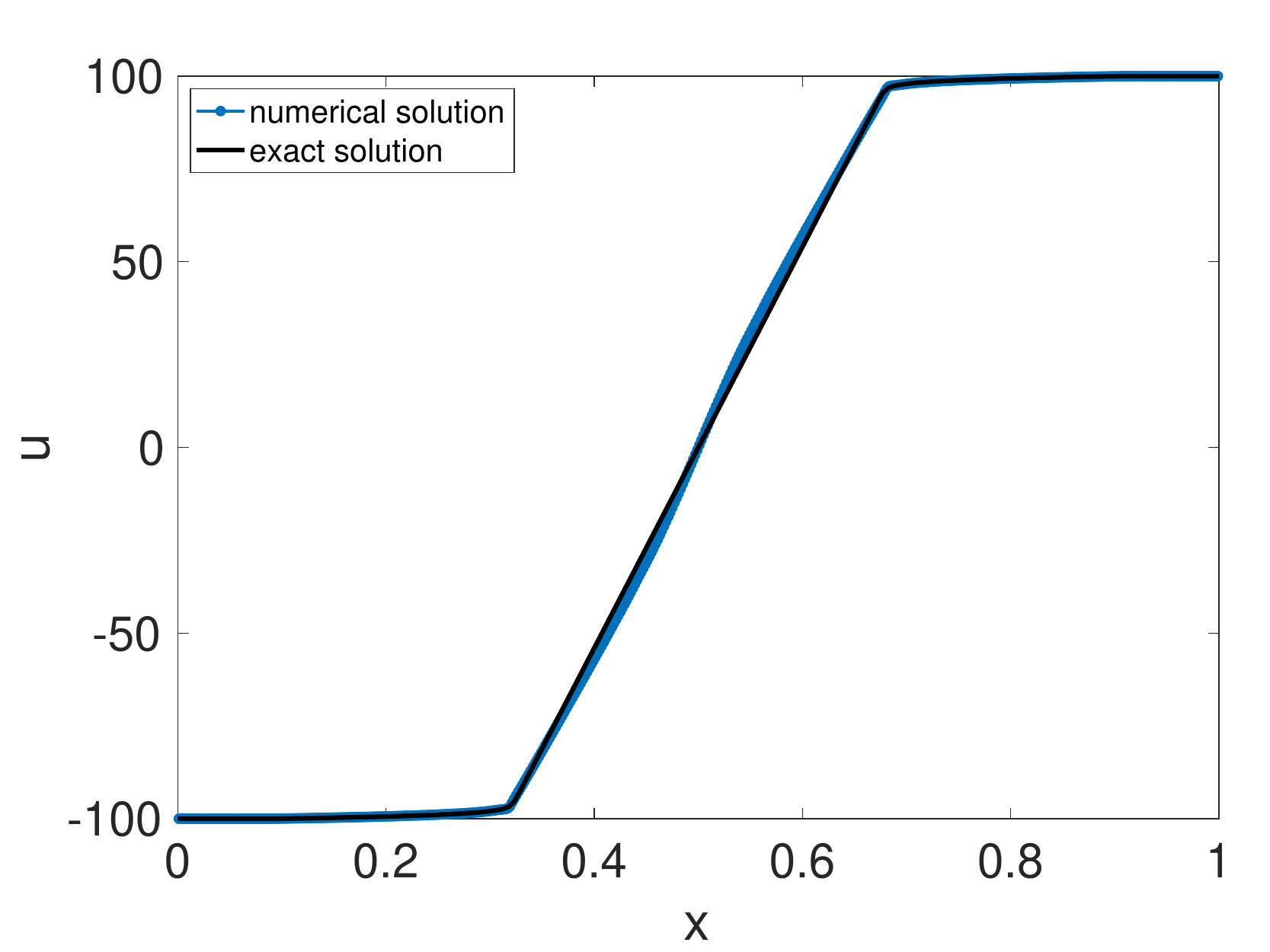}
\label{fig:cavitation-velocity}
}
\caption{Numerical solutions (dots), in comparison with exact ones (black lines), of the cavitation test. The mixture density, pressure, gas volume fraction, and velocity are displayed in subfigures.}
\label{fig:cavitation}
\end{figure}

\vspace{2mm}
\n{\bf Example 4. Cavitation test.} This is the case considered in Ref.~\onlinecite{saurel-simple}. Once again, the mixture of water and air occupies the interval $[0,1]$ and the EOSs of water and air are the same as those used in Example 2. The initial condition is given as follows,
\beqs
(\al_1, \ \rho_1, \ \rho_2, \ u, \ p)=
\left\{
\bga{lr}
(0.99, \ 1000, \ 1, \ -100, \ 10^5), & x<0.5,\\[2mm]
(0.99, \ 1000, \ 1, \ \ \ 100, \ 10^5), & x>0.5.
\eda
\right.
\eeqs
This case demonstrates the ability of the Kapila model to create an interface which does not exist initially.
Perform the numerical computation with $500$ uniform cells to the time $T=1.85$.
The numerical results are displayed in Fig. \ref{fig:cavitation} in comparison with exact ones, where a good agreement is observed.

\begin{figure}[!htb]
\centering
\includegraphics[width=.7\linewidth]{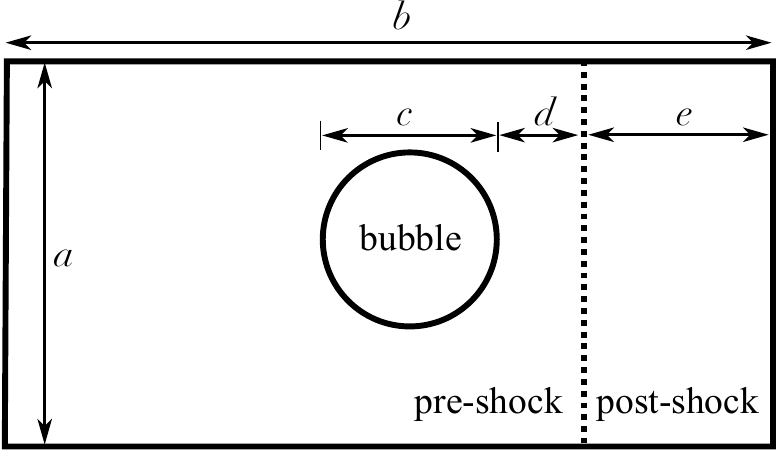}
\caption{The computational domain of the shock-bubble interaction. The dashed line indicates the initial position of the left-moving shock.}
\label{fig:diagram}
\end{figure}

\vspace{2mm}
\n{\bf Example 5. Interaction of a shock in air with a cylindrical helium bubble.} This is a 2D test case \cite{Tryggvason-2009,Quirk-bubble} with experimental data presented in Ref.~\onlinecite{haas}. 
The computational domain is divided into regions shown in Fig. \ref{fig:diagram}, with $a=89$ cm, $b=329.3$ cm, $c=50$ cm, $d=25$ cm, and $e=100$ cm. The dashed line indicates the initial position of the left-moving shock.
The cylindrical bubble is filled with helium while the rest part of the computational domain is filled with air.
Both air and helium are considered as ideal gases with $\gm_1=1.4$ for air and $\gm_2=1.648$ for helium.
The pre-shock air and the bubble remain at rest initially and the pre-shock pressure is $1$ atm. The initial density of the bubble is $0.1819~\text{kg}/\text{m}^3$ and the initial density of the pre-shock air is $1~\text{kg}/\text{m}^3$. The post-shock state is determined by the shock Mach number $\text{Ma}=1.22$.

\begin{figure}[!htb]
\centering
\subfigure[$62\mu s$]{
\includegraphics[width=.44\linewidth]{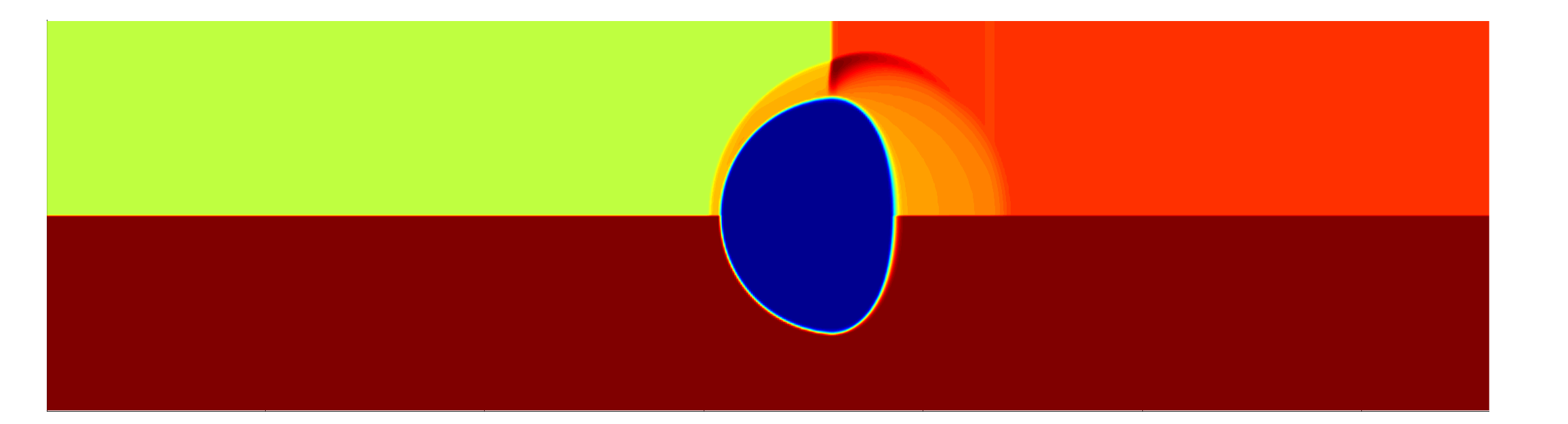}
\label{fig:he-62}
}
\subfigure[$245\mu s$]{
\includegraphics[width=.44\linewidth]{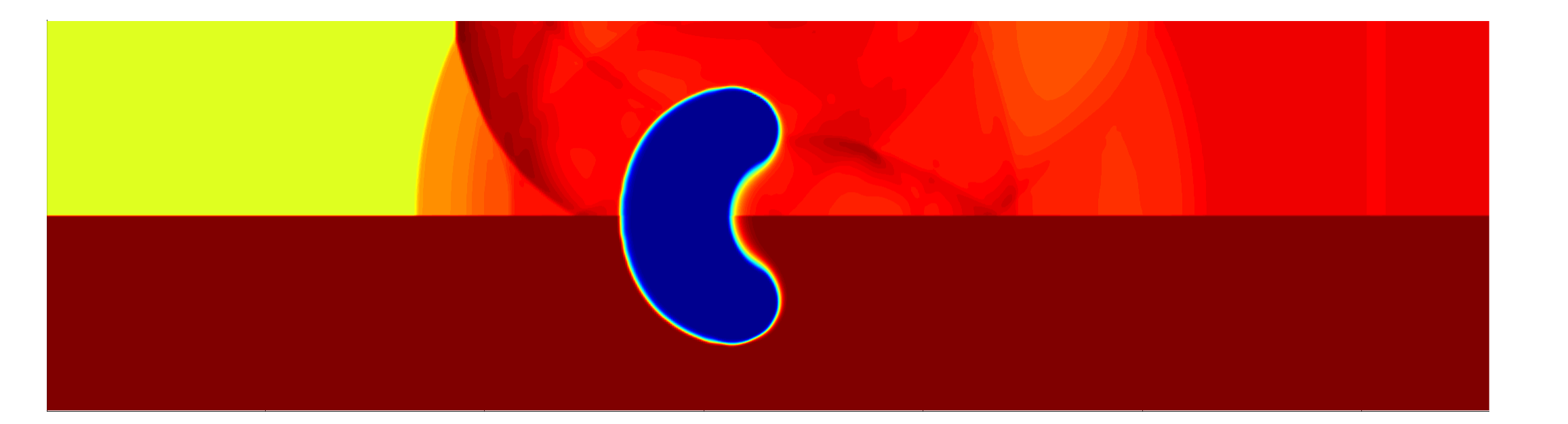}
\label{fig:he-245}
}
\subfigure[$427\mu s$]{
\includegraphics[width=.44\linewidth]{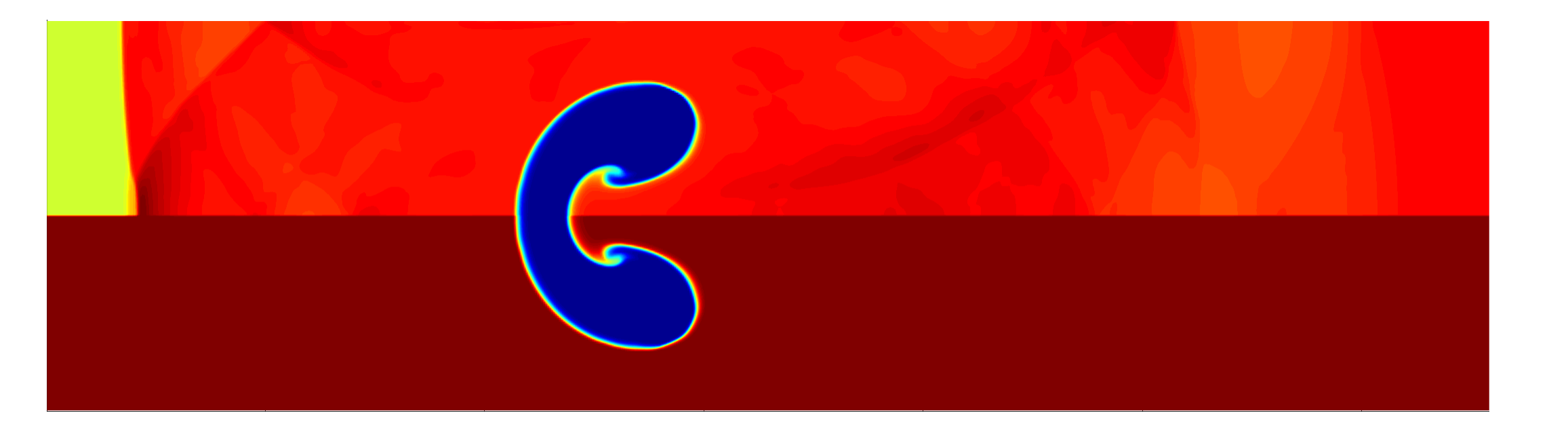}
\label{fig:he-427}
}
\subfigure[$983\mu s$]{
\includegraphics[width=.44\linewidth]{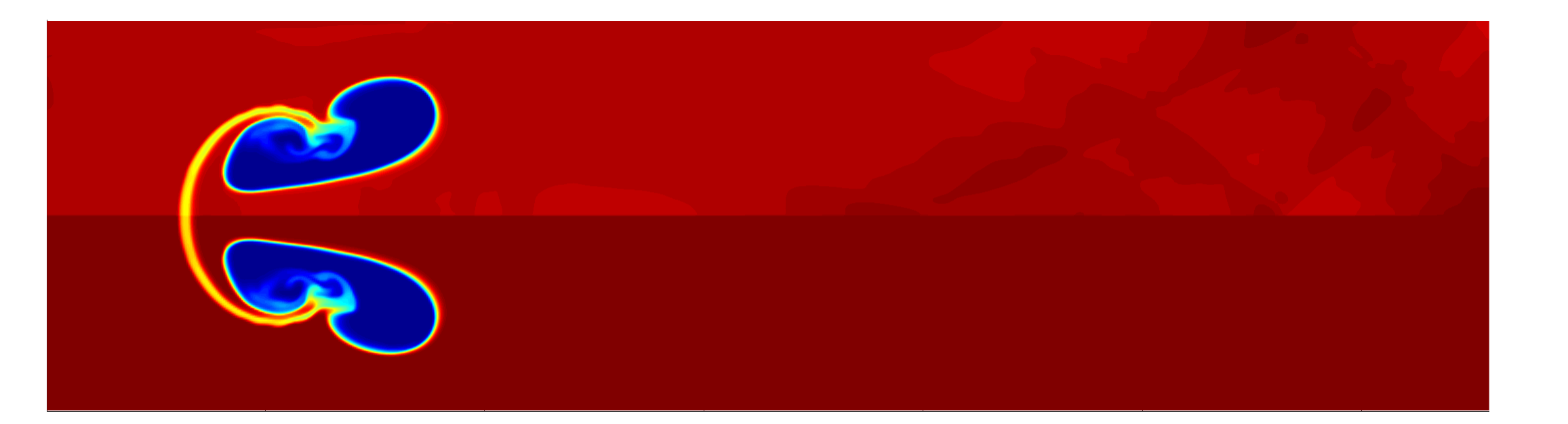}
\label{fig:he-983}
}
\caption{Interaction of a shock in air with a helium cylinder. The mixture density (upper half of each figure) and the air volume fraction (lower half of each figure) are shown at $62\mu s$, $245\mu s$, $427\mu s$, and $983\mu s$ after the shock impact the cylinder.}
\label{fig:helium}
\end{figure}

The numerical computation is performed using $740\times 100$ uniform cells in the upper half of the computational domain due to the symmetry of the flow field. Pseudo-color plots of the mixture density and the air volume fraction are shown in Fig \ref{fig:helium} at $62\mu s$, $245\mu s$, $427\mu s$, and $983\mu s$ after the shock impacts the cylinder. A good agreement with the experimental results in Ref.~\onlinecite{haas} is observed.

\vspace{2mm}
\n{\bf Example 6. Interaction of a shock in water with a cylindrical air bubble.} This is another widely considered two-dimensional case \cite{Tryggvason-2009, LiuTG}.  Here we use the set-up in Ref.~\onlinecite{LiuTG}.
The computational domain is divided into regions demonstrated in Fig. \ref{fig:diagram}, with $a=12$, $b=12$, $c=6$, $d=2.4$, and $e=0.6$.
The cylindrical bubble is filled with air while the rest part of the computational domain is filled with water.
The pre-shock water and the bubble remain at rest initially and the pre-shock pressure is $1$. The initial density of the pre-shock water is $1$ and that of the bubble is $0.0012$. The post-shock state is determined by the shock Mach number $\text{Ma}=1.72$.
Both water and air follow the stiffened gas EOS \eqref{eq:eos-sg}. The parameters for water are taken to be $\gm_1=4.4$ and $\pi_1=6000$, while those for air are taken to be $\gm_2=1.4$ and $\pi_2=0$.

\begin{figure}[!htb]
\centering
\includegraphics[width=.31\linewidth]{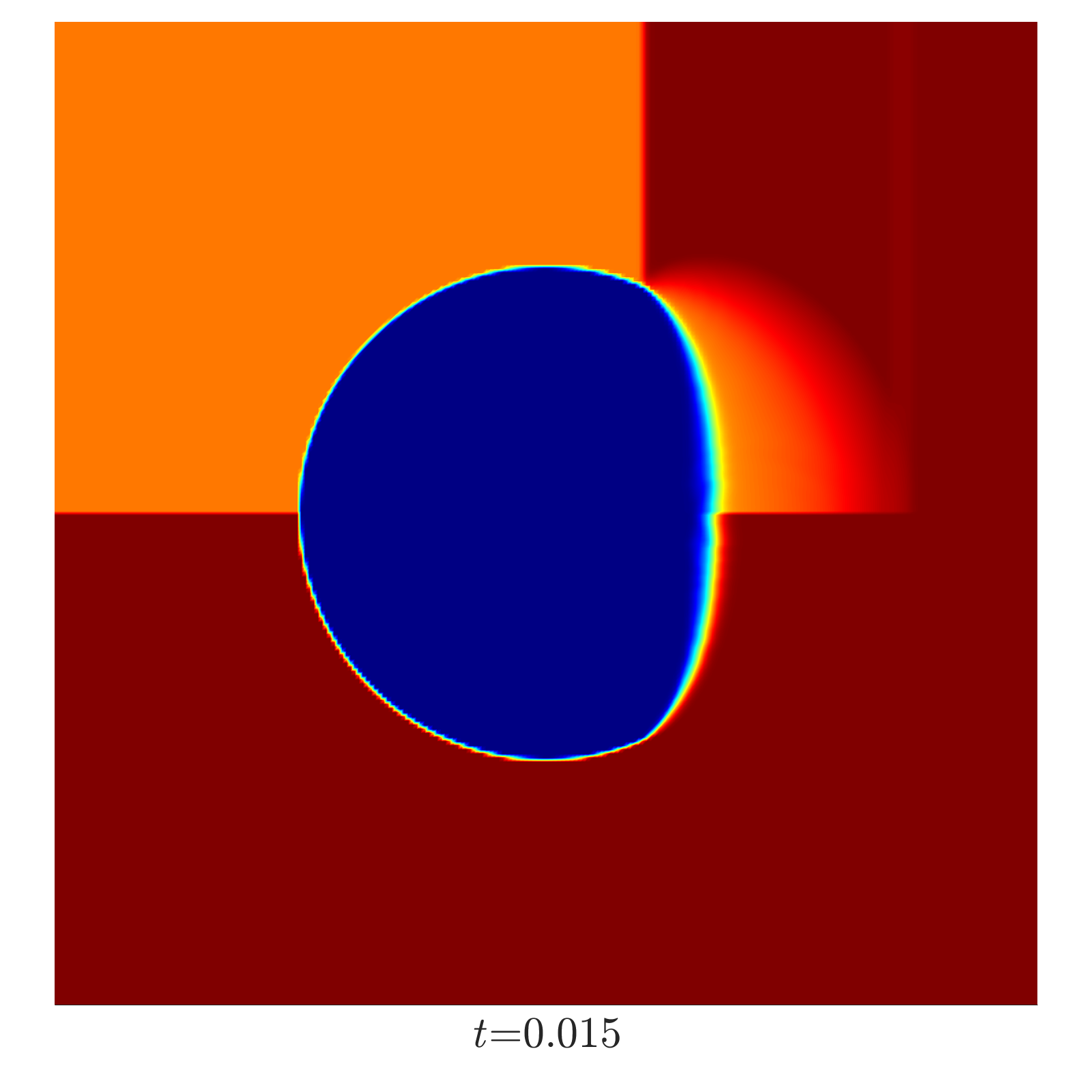}
\includegraphics[width=.31\linewidth]{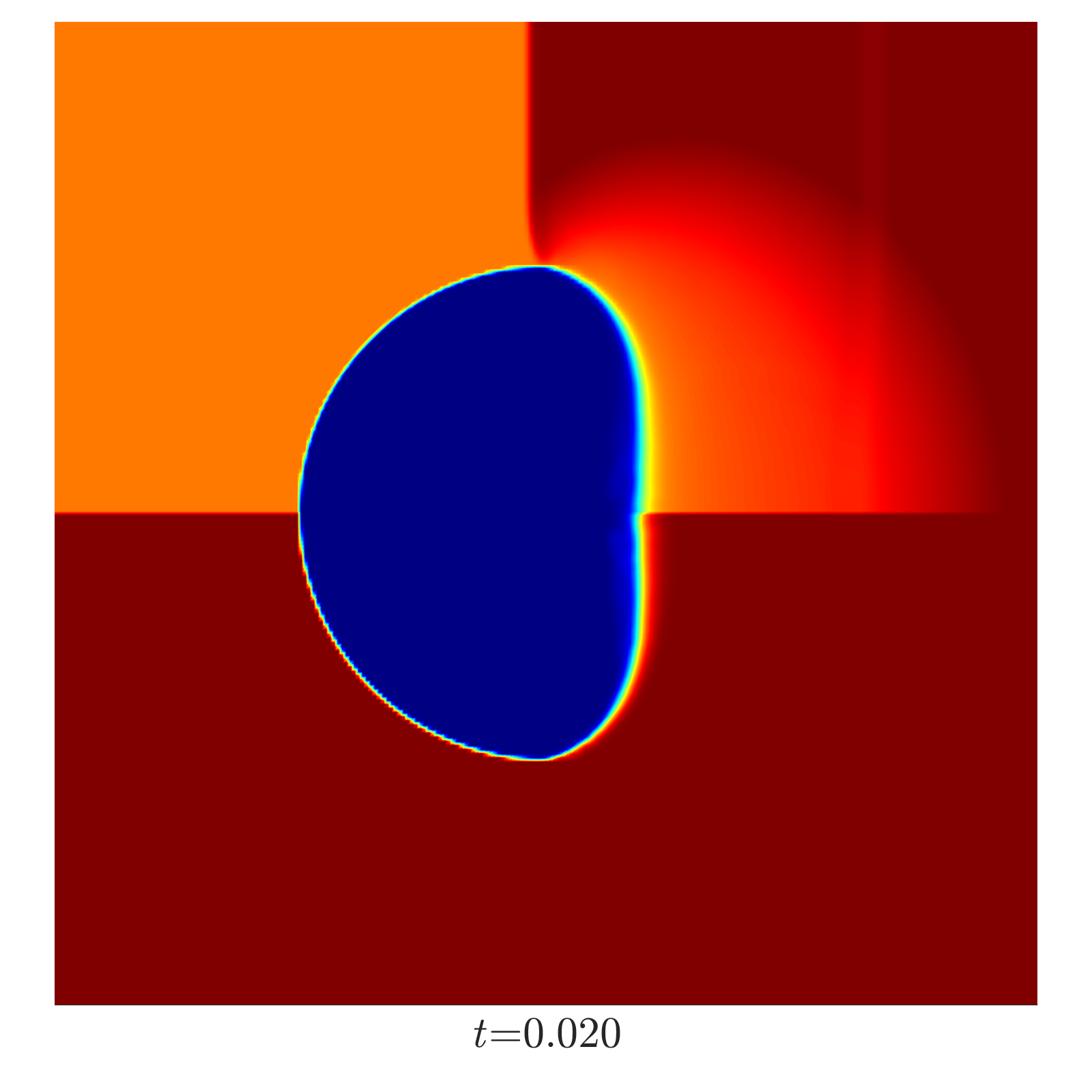}
\includegraphics[width=.31\linewidth]{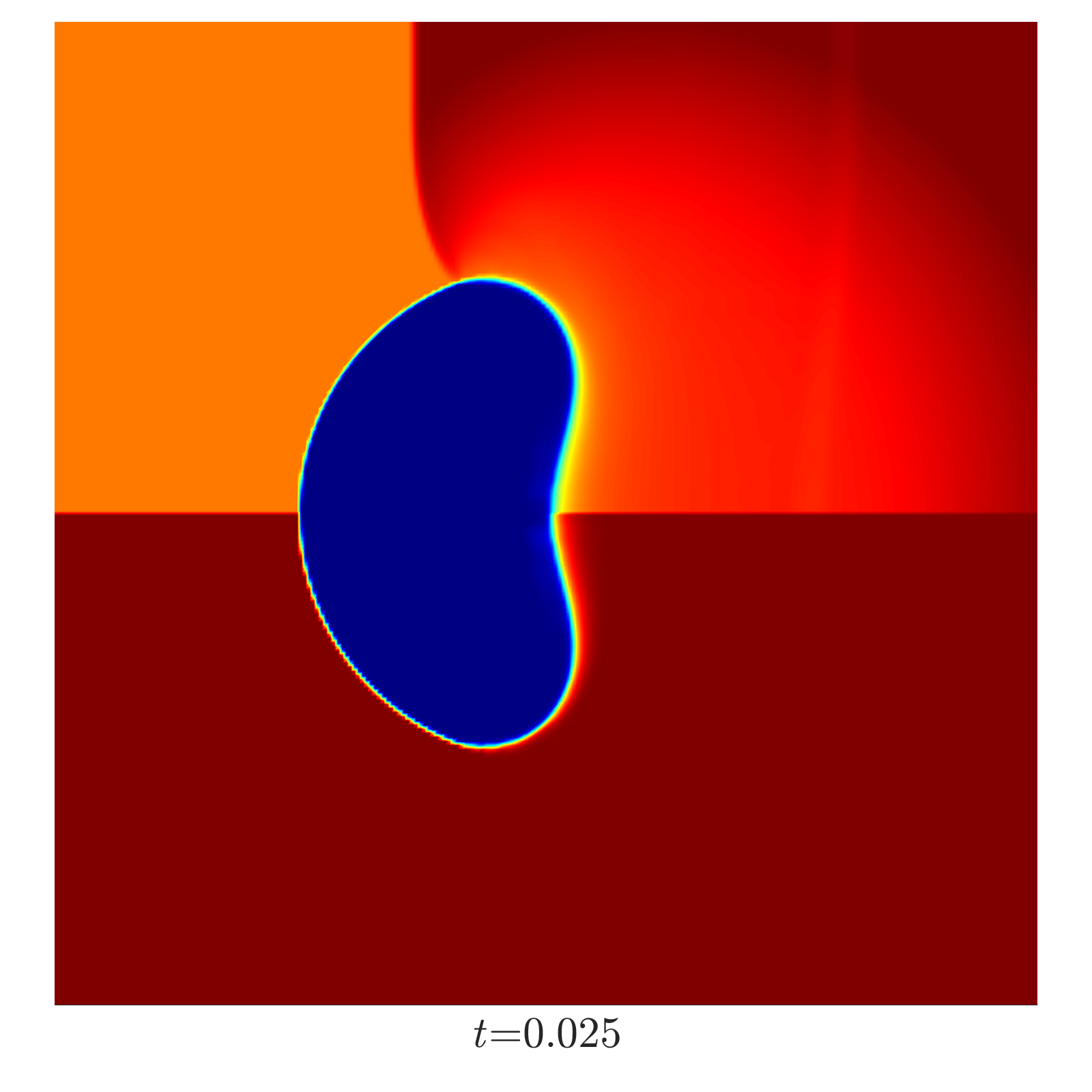}\\[3mm]
\includegraphics[width=.31\linewidth]{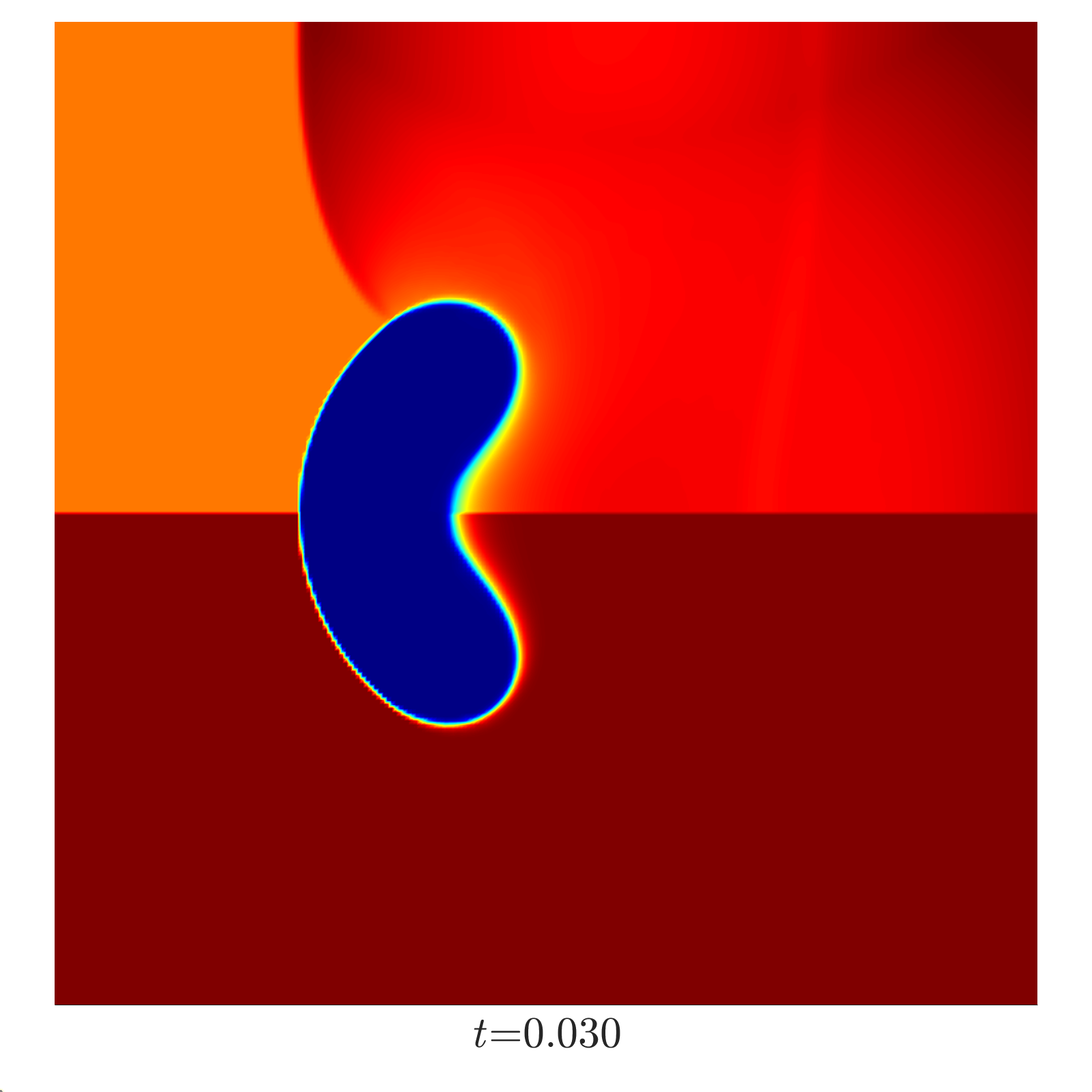}
\includegraphics[width=.31\linewidth]{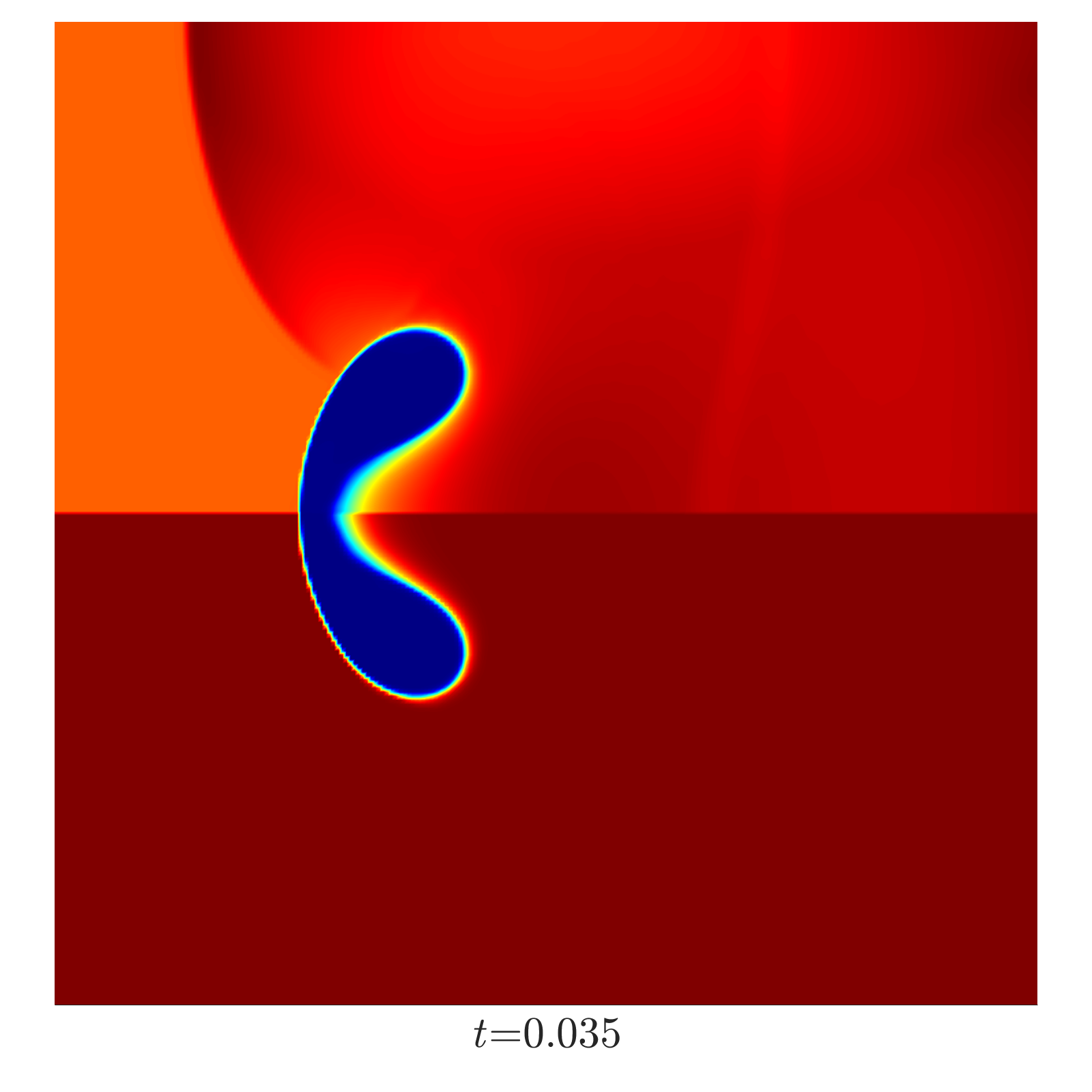}
\includegraphics[width=.31\linewidth]{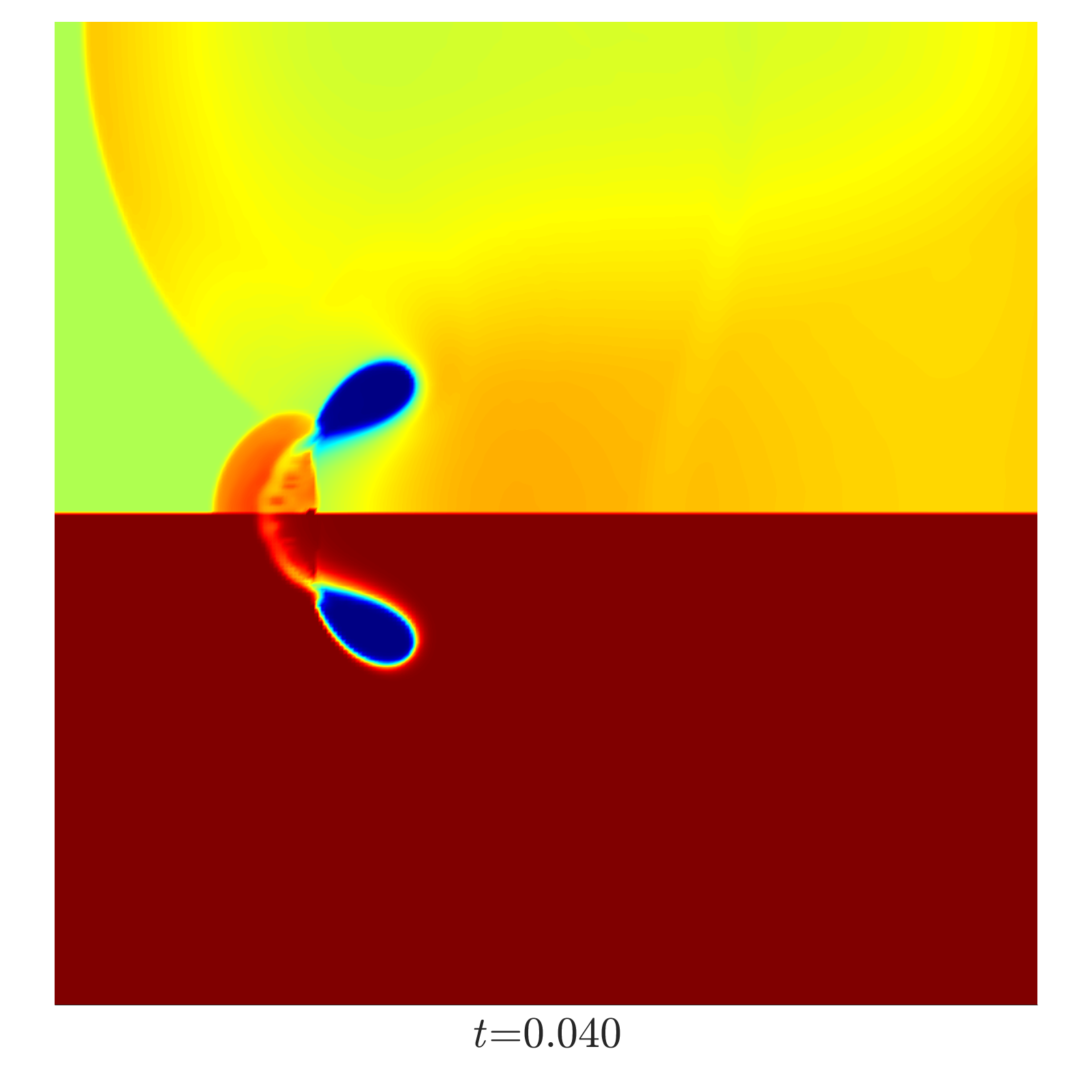}
\caption{Interaction of a shock in water with a cylindrical air bubble. The mixture density (upper half of each figure) and the water volume fraction (lower half of each figure) are shown at computational times $t=0.015$, $t=0.020$, $t=0.025$, $t=0.030$, $t=0.035$, and $t=0.040$, respectively.}
\label{fig:water-air-cylinder}
\end{figure}

By considering the symmetry of the flow field, we perform the numerical computation in the upper half of the computational domain with $240\times 120$ uniform cells. Pseudo-color plots of the mixture density and the water volume fraction are displayed in Fig \ref{fig:water-air-cylinder} at computational times $t=0.015$, $t=0.020$, $t=0.025$, $t=0.030$, $t=0.035$, and $t=0.040$, respectively. A good agreement with the numerical results in \cite{LiuTG} is observed.

\section{Conclusions}\label{sec:discussions}
The Kapila model, derived from the seven-equation model by asymptotic analysis, is widely used in numerical simulations of compressible multiphase flows involving both fluid mixtures and pure fluids \cite{kapila-puckett,murrone,kapila-netherland,kapila-tian,saurel-projection-I, saurel-projection-II}.
By adopting the GRP solver, the present paper contributes to the development of a second-order accurate finite volume scheme for the Kapila model, within the GRP framework \cite{tsc-li,10.1063/1.5113986}.
Given the challenges posed by the stiffness of the Kapila model, the proposed scheme enhances the robustness through two key improvements: a semi-implicit time discretization to the stiff source term and the use of GRP-based numerical fluxes.
The robustness of the current scheme is demonstrated by several challenging test cases involving either interfaces exhibiting large pressure and density ratios or multiphase fluid mixtures. A nonlinear multiphase flow test case with a smooth solution is proposed to verify the numerical accuracy.

It is worth mentioning that in the presence of strong nonlinear waves, the nonlinear GRP solver is imperative since it effectively resolves essential thermodynamical information regarding the evolution of the flow field, as demonstrated in the shock tube problem of Example 2 in Section \ref{sec:numer}.
To keep the conciseness of this paper, we left this topic to a separate study Ref.~\onlinecite{grp-kapila}.


\appendix
\section{Space data reconstruction}\label{sec:app}
To complete the finite volume scheme, the linear data reconstruction employed in Refs.~\onlinecite{benartzi-01,grp-06} is adopted. Here we apply the linear reconstruction to primitive variables $\bW=[\zeta_1,\rho,u,v,p,\al_1]^\top$.

In two space dimensions, $\bW$ is linearly reconstructed in cell $\Omega_{ij}$ as\beql{eq:recon-2d}
\bW^{n+1}_{i,j}(\bx)=\widetilde\bW_{i,j}^{n+1}+(\nabla\bW)_{i,j}^{n+1}\cdot(\bx-\bx_{i,j}).
\eeq
The cell centered value $\widetilde\bW_{i,j}^{n+1}$ is computed using cell average values $\overline\bU_{i,j}^{n+1}$ obtained in \eqref{eq:scheme}. The approximation to the slope is
\beql{eq:slope-def}
\left(\dfr{\pt\bW}{\pt x}\right)_{i,j}^{n+1}=\minmod
\left(
\dfr{\widetilde\bW_{i,j}^{n+1}-\widetilde\bW_{i-1,j}^{n+1}}{\Dx},
\kappa\dfr{\widehat\bW_{i+\frac 12,j}^{n+1}-\widehat\bW_{i-\frac 12,j}^{n+1}}{\Dx},
\dfr{\widetilde\bW_{i+1,j}^{n+1}-\widetilde\bW_{i,j}^{n+1}}{\Dx}
\right),
\eeq
where the cell interface values are computed as
\beql{eq:interface-inte-numer}
\widehat\bW_{i+\frac 12,j}^{n+1}=\d\sum_{g}\omega_g\widehat\bW_{i+\frac 12,j}^{n+1,g}, \ \
\widehat\bW_{i+\frac 12,j}^{n+1} = \bW_{i+\frac 12,j}^{n,*,g} + \left(\dfr{\pt\bW}{\pt t}\right)_{i+\frac 12,j}^{n,*,g}.
\eeq
The coefficient $\kappa\in[0,2)$ is a user-tuned parameter.
The minmod function is defined as
\beql{eq:minmod-def}
\minmod(a,b,c)=
\left\{
\bga{ll}
\text{sign}(a) \min(|a|,|b|,|c|), & \text{if }\text{sign}(a)=\text{sign}(b)=\text{sign}(c), \\
0, & \text{otherwise}.
\eda
\right.
\eeq
The other component $(\pt\bW/\pt y)_{i,j}^{n+1}$ of the gradient can be computed similarly. 
By transforming the reconstructed piecewise linear function of $\bW$ into those of $\bU$, the initial value to be used in the IVP \eqref{eq:grp-def} is obtained.

\begin{acknowledgments}
This research is supported by NSFC (Nos. 11901045, 12031001), and Hubei Ed. Dept. Sci. $\&$ Tech. Res. Project (No. D20232901).
\end{acknowledgments}
\vspace{2mm}

\bibliographystyle{siam}
\bibliography{KapilaFV.bbl}

\begin{thebibliography}{10}

\bibitem{star-mre}
{\sc B.~Albertazzi, P.~Mabey, T.~Michel, G.~Rigon, J.~R. Marqu\`es, S.~Pikuz,
  S.~Ryazantsev, E.~Falize, L.~V.~B. Som, J.~Meinecke, N.~Ozaki, G.~Gregori,
  and M.~Koenig}, {\em Triggering star formation: {E}xperimental compression of
  a foam ball induced by {T}aylor–{S}edov blast waves}, Matter Radiat.
  Extremes, 7 (2022), p.~036902.

\bibitem{bn-saurel-warnecke}
{\sc N.~Andrianov, R.~Saurel, and G.~Warnecke}, {\em A simple method for
  compressible multiphase mixtures and interfaces}, Int. J. Numer. Methods
  Fluids, 41 (2003), pp.~109--131.

\bibitem{bn-warnecke}
{\sc N.~Andrianov and G.~Warnecke}, {\em The {R}iemann problem for the
  {B}aer-{N}unziato two-phase flow model}, J. Comput, Phys., 212 (2004),
  pp.~434--464.

\bibitem{b-n}
{\sc M.~R. Baer and J.~W. Nunziato}, {\em Two-phase modeling of
  deflagration-to-detonation transition in granular materials: {A} critical
  examination of modeling issues}, Int. J. Multiphase Flow, 12 (1986),
  pp.~861--889.

\bibitem{benartzi-84}
{\sc M.~Ben-Artzi and J.~Falcovitz}, {\em A second-order {G}odunov-type scheme
  for compressible fluid dynamics}, J. Comput. Phys., 55 (1984), pp.~1--32.

\bibitem{benartzi-01}
{\sc M.~Ben-Artzi and J.~Falcovitz}, {\em Generalized Riemann Problems in
  Computational Fluid Dynamics}, Cambridge University Press, Cambridge, 2003.

\bibitem{grp-06}
{\sc M.~Ben-Artzi, J.~Li, and G.~Warnecke}, {\em Direct {E}ulerian {GRP} scheme
  for compressible fluid flows}, J. Comput, Phys., 31 (2006), pp.~335--362.

\bibitem{Fedkiw}
{\sc R.~Caiden, R.~P. Fedkiw, and C.~Anderson}, {\em A numerical method for
  two-phase flow consisting of separate compressible and incompressible
  regions}, J. Comput. Phys., 166 (2001), pp.~1--27.

\bibitem{detonation-saurel-jcp}
{\sc A.~Chinnayya, E.~Daniel, and R.~Saurel}, {\em Modelling detonation waves
  in heterogeneous energetic materials}, J. Comput. Phys., 196 (2004),
  pp.~769--774.

\bibitem{grp-kapila}
{\sc Z.~Du}, {\em Nonlinear generalized {R}iemann problem solver for the
  {K}apila model of compressible multiphase flows}.
\newblock Paper in preparation, 2025.

\bibitem{du-hweno}
{\sc Z.~Du and J.~Li}, {\em A {H}ermite {WENO} reconstruction for fourth order
  temporal accurate schemes based on the {GRP} solver for hyperbolic
  conservation laws}, J. Comput, Phys., 355 (2018), pp.~3045--3069.

\bibitem{Toro-ADER-pnpm}
{\sc M.~Dumbser, D.~S. Balsara, E.~F. Toro, and C.-D. Munz}, {\em A unified
  framework for the construction of one-step finite volume and discontinuous
  {G}alerkin schemes on unstructured meshes}, J. Comput. Phys., 227 (2008),
  pp.~8209--8253.

\bibitem{Osher-2018}
{\sc R.~F. F.~Gibou and S.~Osher}, {\em A review of level-set methods and some
  recent applications}, J. Comput. Phys., 353 (2018), pp.~82--109.

\bibitem{haas}
{\sc J.-F. Haas and B.~Sturtevant}, {\em Interaction of weak shock waves with
  cylindrical and spherical gas inhomogeneities}, J. Fluid Mech., 181 (1987),
  pp.~41--76.

\bibitem{LiXL-2015}
{\sc Y.~Hu, Q.~Shi, V.~F. De~Almeida, and X.~Li}, {\em Numerical simulation of
  phase transition problems with explicit interface tracking}, Chem. Eng. Sci.,
  128 (2015), pp.~92--108.

\bibitem{gks-hweno-0}
{\sc X.~Ji, L.~Pan, W.~Shyy, and K.~Xu}, {\em A compact fourth-order
  gas-kinetic scheme for the {E}uler and {N}avier-{S}tokes equations}, J.
  Comput, Phys., 372 (2018), pp.~446--472.

\bibitem{gks-hweno}
{\sc X.~Ji, F.~Zhao, W.~Shyy, and K.~Xu}, {\em A {HWENO} reconstruction based
  high-order compact gas-kinetic scheme on unstructured mesh}, J. Comput,
  Phys., 410 (2020), p.~109367.

\bibitem{kamaya-1996}
{\sc H.~Kamaya}, {\em On the origin of small-scale structure in
  self-gravitating two-phase gas}, Astrophys. J., 465 (1996), pp.~769--774.

\bibitem{kapila}
{\sc A.~K. Kapila, R.~Menikoff, J.~B. Bdzil, and S.~F. Son}, {\em Two-phase
  modeling of deflagration-to-detonation transition in granular materials:
  {R}educed equations}, Phys. Fluids, 13 (2001), pp.~3002--3024.

\bibitem{kapila-structure}
{\sc A.~K. Kapila, S.~F. Son, J.~B. Bdzil, R.~Menikoff, and D.~S. Stewart},
  {\em Two-phase modeling of {DDT}: {S}tructure of the velocity-relaxation
  zone}, Phys. Fluids, 9 (1997), pp.~3885--3897.

\bibitem{bn-karni}
{\sc S.~Karni and G.~Hern\'andez-Due\~nas}, {\em A hybrid algorithm for the
  {B}aer-{N}unziato model using the {R}iemann invariants}, J. Sci. Comput., 45
  (2010), pp.~382--403.

\bibitem{kapila-netherland}
{\sc J.~J. Kreeft and B.~Koren}, {\em A new formulation of {K}apila’s
  five-equation model for compressible two-fluid flow, and its numerical
  treatment}, J. Comput, Phys., 229 (2010), pp.~6220--6242.

\bibitem{L-W}
{\sc P.~D. Lax and B.~Wendroff}, {\em Systems of conservation laws}, Comm. Pure
  Appl. Math., 13 (1960), pp.~217--237.

\bibitem{cavitation-saurel-jcp}
{\sc O.~Le~Metayer, J.~Massoni, and R.~Saurel}, {\em Modelling evaporation
  fronts with reactive {R}iemann solvers}, J. Comput. Phys., 205 (2005),
  pp.~567--610.

\bibitem{Lei-Li-2018-AMM}
{\sc X.~Lei and J.~Li}, {\em Transversal effects and genuine
  multidimensionality of high order numerical schemes for compressible fluid
  flows}, Appl. Math. Mech. -Engl. Ed., 39 (2018), pp.~1--12.

\bibitem{bn-lei}
{\sc X.~Lei and J.~Li}, {\em A staggered-projection {G}odunov-type method for
  the {B}aer-{N}unziato two-phase model}, J. Comput. Phys., 437 (2021),
  p.~110312.

\bibitem{tsc-li}
{\sc J.~Li}, {\em Two-stage fourth order: temporal-spatial coupling in
  computational fluid dynamics ({CFD})}, Adv. Aerodyn, 1 (2019), pp.~1--36.

\bibitem{10.1063/1.5113986}
{\sc J.~Li, X.~Lei, Z.~Du, and Y.~Wang}, {\em High order temporal-spatially
  coupled schemes for compressible multi-fluid flows}, AIP Conference
  Proceedings, 2116 (2019), p.~030002.

\bibitem{combustion}
{\sc M.~A. Liberman}, {\em Combustion Physics: Flames, Detonations, Explosions,
  Astrophysical Combustion and Inertial Confinement Fusion}, Springer Nature,
  Switzerland, 2021.

\bibitem{menshov-1991}
{\sc I.~S. Men'shov}, {\em The generalized problemm of breakup of an arbitrary
  discontinuity}, J. Appl. Maths Mechs., 55 (1991), pp.~67--74.

\bibitem{kapila-puckett}
{\sc G.~H. Miller and E.~G. Puckett}, {\em A high-order {G}odunov method for
  multiple condensed phases}, J. Comput, Phys., 128 (1996), pp.~134--164.

\bibitem{Toro-ADER-blood}
{\sc G.~I. Montecinos, A.~Santac\'a, M.~Celant, L.~O. M\"uller, and E.~F.
  Toro}, {\em A unified framework for the construction of one-step finite
  volume and discontinuous {G}alerkin schemes on unstructured meshes}, Comput.
  \& Fluids, 248 (2022), p.~105685.

\bibitem{murrone}
{\sc A.~Murrone and H.~Guillard}, {\em A five equation reduced model for
  compressible two phase flow problems}, J. Comput, Phys., 202 (2005),
  pp.~664--698.

\bibitem{saurel-projection-II}
{\sc F.~Petitpas, E.~Franquet, R.~Saurel, and O.~Le~Metayer}, {\em A
  relaxation-projection method for compressible flows. {P}art {II}:
  {A}rtificial heat exchanges for multiphase shocks}, J. Comput. Phys., 225
  (2007), pp.~2214--2248.

\bibitem{saurel-ijmf}
{\sc F.~Petitpas, J.~Massoni, R.~Saurel, E.~Lapebie, and L.~Munier}, {\em
  Diffuse interface model for high speed cavitating underwater systems}, Int.
  J. Multiphase Flows, 35 (2009), pp.~747--759.

\bibitem{detonation-saurel-sw}
{\sc F.~Petitpas, R.~Saurel, E.~Franquet, and A.~Chinnayya}, {\em Modelling
  detonation waves in condensed materials: multiphase cj conditions and
  multidimensional computations}, Shock Waves, 19 (2009), pp.~377--401.

\bibitem{Quirk-bubble}
{\sc J.~J. Quirk and S.~Karni}, {\em On the dynamics of a shock-bubble
  interaction}, J . Fluid Mech., 318 (1996), pp.~129--163.

\bibitem{bn-abgrall}
{\sc R.~Saurel and R.~Abgrall}, {\em A multiphase {G}odunov method for
  compressible multifluid and multiphase flows}, J. Comput, Phys., 150 (1999),
  pp.~425--467.

\bibitem{saurel-projection-I}
{\sc R.~Saurel, E.~Franquet, E.~Daniel, and O.~Le~Metayer}, {\em A
  relaxation-projection method for compressible flows. {P}art {I}: {T}he
  numerical equation of state for {E}uler equations}, J. Comput. Phys., 223
  (2007), pp.~822--845.

\bibitem{detonation-saurel-jfm}
{\sc R.~Saurel and O.~Le~Metayer}, {\em A multiphase model for interfaces,
  shocks, detonation waves and cavitation}, J. Comput. Phys., 431 (2001),
  pp.~239--271.

\bibitem{saurel-shock-jump}
{\sc R.~Saurel, O.~Le~Metayer, J.~Massoni, and S.~Gavrilyuk}, {\em Shock jump
  relations for multiphase mixtures with stiff mechanical relaxation}, Shock
  Waves, 16 (2007), pp.~209--232.

\bibitem{review-saurel}
{\sc R.~Saurel and C.~Pantano}, {\em Diffuse-interface capturing methods for
  compressible two-phase flows}, Annu. Rev. Fluid Mech., 50 (2018),
  pp.~105--130.

\bibitem{saurel-jfm}
{\sc R.~Saurel, F.~Petitpas, and R.~Abgrall}, {\em Modelling phase transition
  in metastable liquids: application to cavitating and flashing flows}, J.
  Fluid Mech., 607 (2008), pp.~313--350.

\bibitem{saurel-simple}
{\sc R.~Saurel, F.~Petitpas, and R.~A. Berry}, {\em Simple and efficient
  relaxation methods for interfaces separating compressible fluids, cavitating
  flows and shocks in multiphase mixtures}, J. Comput, Phys., 228 (2009),
  pp.~1678--1712.

\bibitem{bn-kapila}
{\sc D.~W. Schwendeman, C.~W. Wahle, and A.~K. Kapila}, {\em The {R}iemann
  problem and a high-resolution {G}odunov method for a model of compressible
  two-phase flow}, J. Comput. Phys., 212 (2006), pp.~490--526.

\bibitem{Tryggvason-2009}
{\sc H.~Terashima and G.~Tryggvason}, {\em A front-tracking/ghost-fluid method
  for fluid interfaces in compressible flows}, J. Coumpt. Phys., 228 (2009),
  pp.~4012--4037.

\bibitem{kapila-tian}
{\sc B.~Tian and L.~Li}, {\em A five-equation model based global {ALE} method
  for compressible multifluid and multiphase flows}, Comput. \& Fluids, 241
  (2021), pp.~149--172.

\bibitem{Toro-ADER}
{\sc E.~F. Toro and V.~A. Titarev}, {\em Derivative {R}iemann solvers for
  systems of conservation laws and {ADER} methods}, J. Comput. Phys., 212
  (2006), pp.~150--165.

\bibitem{wood}
{\sc A.~Wood}, {\em A Textbook of Sound}, MacMillan, New York, 1930.

\bibitem{LiuTG}
{\sc L.~Xu and T.~Liu}, {\em Explicit interface treatments for compressible
  gas-liquid simulations}, Comput. \& Fluids, 153 (2017), pp.~34--48.

\end{thebibliography}

\section*{Data Availability Statement}
The data that support the findings of this study are available from the corresponding author upon reasonable request.

\end{document}